\documentclass[11pt]{article}
\usepackage{latexsym,amssymb}
\usepackage{amsmath}

\newtheorem{Theorem}{\bf Theorem}[section]
\newtheorem{Lemma}{\bf Lemma}[section]
\newtheorem{Proposition}{\bf Proposition}[section]
\newtheorem{Corollary}{\bf Corollary}[section]
\newtheorem{Remark}{\bf Remark}[section]
\newtheorem{Example}{\bf Example}[section]
\newtheorem{Definition}{\bf Definition}[section]

\newenvironment{theorem}{\begin{Theorem}$\!\!\!$}{\end{Theorem}}
\newenvironment{lemma}{\begin{Lemma}$\!\!\!$}{\end{Lemma}}

\newenvironment{corollary}{\begin{Corollary}$\!\!\!$}{\end{Corollary}}
\newenvironment{remark}{\begin{Remark}$\!\!\!$}{\end{Remark}}

\newenvironment{definition}{\begin{Definition}$\!\!\!$}{\end{Definition}}

\numberwithin{equation}{section}

\begin{document}

\title{Heat equation with a nonlinear boundary condition\\
and uniformly local $L^r$ spaces}
\author{Kazuhiro Ishige and Ryuichi Sato\\
Mathematical Institute, Tohoku University\\
Aoba, Sendai 980-8578, Japan
}
\date{}
\maketitle
\begin{abstract}
We establish the local existence and the uniqueness of solutions 
of the heat equation with a nonlinear boundary condition 
for the initial data in uniformly local $L^r$ spaces.  
Furthermore, we study the sharp lower estimates of the blow-up time 
of the solutions with the initial data $\lambda\psi$ as $\lambda\to 0$ or $\lambda\to\infty$ 
and the lower blow-up estimates of the solutions. 
\end{abstract}
\section{Introduction}
This paper is concerned with the heat equation with a nonlinear boundary condition, 
\begin{equation}
\left\{
\begin{array}{ll}
\partial_t u=\Delta u, & x\in\Omega,\,t>0,\vspace{3pt}\\
\nabla u\cdot\nu(x)=|u|^{p-1}u,\qquad &x\in\partial\Omega,\,\,t>0,\vspace{3pt}\\
u(x,0)=\varphi(x), & x\in\Omega,
\end{array}
\right.
\label{eq:1.1}
\end{equation}
where $N\ge 1$, $p>1$,
$\Omega$ is a smooth domain in ${\bf R}^N$, $\partial_t=\partial/\partial t$ and 
$\nu=\nu(x)$ is the outer unit normal vector to $\partial\Omega$. 
For any $\varphi\in BUC(\Omega)$, 
problem~\eqref{eq:1.1} has a unique solution 
$$
u\in C^{2,1}(\Omega\times(0,T])\,\cap\,C^{1,0}(\overline{\Omega}\times(0,T])\,\cap\, BUC(\overline{\Omega}\times[0,T])
$$
for some $T>0$ and 
the maximal existence time $T(\varphi)$ of the solution can be defined. 
If $T(\varphi)<\infty$, then 
$$
\limsup_{t\to\,T(\varphi)}\|u(t)\|_{L^\infty(\Omega)}=\infty 
$$
and we call $T(\varphi)$ the blow-up time of the solution $u$. 

Problem~\eqref{eq:1.1} has been studied in many papers
from various points  of view 
(see e.g. \cite{CF01}--\cite{DB01}, \cite{FR}--\cite{GL}, \cite{GH}--\cite{HM}, \cite{IK}, \cite{K}, \cite{QS2} and references therein) 
while there are few results related to the dependence of the blow-up time on the initial function 
even in the case $\Omega={\bf R}^N_+$. 
We remark that 
the blow-up time for problem~\eqref{eq:1.1} cannot be chosen uniform for all initial functions
lying in a bounded set of $L^r({\bf R}^N_+)$ with $1\le r\le N(p-1)$. 
Indeed, similarly to \cite[Remark~15.4~(i)]{QS}, 
for any solution $u$ blowing up at $t=T<\infty$ and $\mu>0$, 
\begin{equation}
\label{eq:1.2}
u_\mu(x,t):=\mu^{1/(p-1)}u(\mu x,\mu^2 t)
\end{equation}
is a solution of \eqref{eq:1.1} blowing up at $t=\mu^{-2}T$ while  
$$
\|u_\mu(0)\|_{L^r({\bf R}^N_+)}
=\mu^{\frac{1}{p-1}-\frac{N}{r}}\|\varphi\|_{L^r({\bf R}^N_+)}
\le\|\varphi\|_{L^r({\bf R}^N_+)}
$$ 
for any $\mu\ge 1$. 
\vspace{5pt}

For $1\le r<\infty$ and $\rho>0$, 
let $L^r_{uloc,\rho}(\Omega)$ be the uniformly local $L^r$ space in $\Omega$ equipped with the norm 
$$
||f||_{r,\rho}:=\sup_{x\in\overline{\Omega}}\,\left(\int_{\Omega\,\cap\,B(x,\rho)}|f(y)|^rdy\right)^{1/r}. 
$$
We denote by ${\cal L}^r_{uloc,\rho}(\Omega)$ 
the completion of bounded uniformly continuous functions in $\Omega$ with respect to the norm $\|\cdot\|_{r,\rho}$, that is, 
$$
{\cal L}^r_{uloc,\rho}(\Omega):=\overline{BUC(\Omega)}^{\,\|\,\cdot\,\|_{r,\rho}}. 
$$
We set $L^\infty_{uloc,\rho}(\Omega)=L^\infty(\Omega)$ and ${\cal L}^\infty_{uloc,\rho}(\Omega)=BUC(\Omega)$.
\vspace{3pt}

In this paper we prove the local existence and the uniqueness of the solutions of problem~\eqref{eq:1.1} 
with initial functions in ${\cal L}^r_{uloc,\rho}(\Omega)$, 
and study the dependence of the blow-up time on the initial functions. 
As an application of the main results of this paper, 
we study the asymptotic behavior of the blow-up time $T(\varphi)$ with $\varphi=\lambda\psi$ 
as $\lambda\to 0$ or $\lambda\to\infty$ and show the validity of our arguments. 
Furthermore, we obtain a lower estimate of the blow-up rate of the solutions 
(see Section~5). 
\vspace{5pt}

Throughout this paper,  
following \cite[Section~1]{QS}, we assume that $\Omega\subset{\bf R}^N$ is a uniformly regular domain of class $C^1$. 
For any $x\in{\bf R}^N$ and $\rho>0$, define
$$
B(x,\rho):=\{y\in{\bf R}^N:|x-y|<\rho\},\,\,\,
\Omega(x,\rho):=\Omega\,\cap\,B(x,\rho),\,\,\,
\partial\Omega(x,\rho):=\partial\Omega\,\cap\,B(x,\rho).
$$
By the trace inequality for $W^{1,1}(\Omega)$-functions 
and the Gagliardo-Nirenberg inequality 
we can find $\rho_*\in(0,\infty]$ with the following properties (see Lemma~\ref{Lemma:2.2}). 
\begin{itemize}
  \item There exists a positive constant $c_1$ such that 
  \begin{equation}
  \label{eq:1.3}
  \int_{\partial\Omega(x,\rho)}|v|\,d\sigma\le c_1\int_{\Omega(x,\rho)}|\nabla v|\,dy
  \end{equation}
  for all $v\in C_0^1(B(x,\rho))$, $x\in\overline{\Omega}$ and $0<\rho<\rho_*$. 
  \item Let $1\le \alpha$, $\beta\le\infty$ and $\sigma\in[0,1]$ be such that 
  \begin{equation}
  \label{eq:1.4}
  \frac{1}{\alpha}=\sigma\left(\frac{1}{2}-\frac{1}{N}\right)+(1-\sigma)\frac{1}{\beta}. 
  \end{equation}
  Assume, if $N\ge 2$, that $\alpha\not=\infty$ or $N\not=2$. 
  Then there exists a constant $c_2$ such that 
  \begin{equation}
  \label{eq:1.5}
  \|v\|_{L^\alpha(\Omega(x,\rho))}\le c_2\|v\|_{L^\beta(\Omega(x,\rho))}^{1-\sigma}\|\nabla v\|_{L^2(\Omega(x,\rho))}^\sigma
  \end{equation}
  for all $v\in C^1_0(B(x,\rho))$, $x\in\overline{\Omega}$ and $0<\rho<\rho_*$. 
\end{itemize} 
We remark that, in the case 
$$
\Omega=\{(x',x_N)\in{\bf R}^N\,:\,x_N>\Phi(x')\}, 
$$
where $N\ge 2$ and $\Phi\in C^1({\bf R}^{N-1})$ with $\|\nabla\Phi\|_{L^\infty({\bf R}^{N-1})}<\infty$, 
\eqref{eq:1.3} and \eqref{eq:1.5} hold with $\rho_*=\infty$
(see Lemma~\ref{Lemma:2.2}).  
Inequalities~\eqref{eq:1.3} and \eqref{eq:1.5} 
are used to treat the nonlinear boundary condition.  
\vspace{5pt}

Next we state the definition of the solution of \eqref{eq:1.1}. 
\begin{definition}
\label{Definition:1.1}
Let $0<T\le\infty$ and $1\le r<\infty$. 
Let $u$ be a continuous function in $\overline\Omega\times(0,T]$. 
We say that 
$u$ is a $L^r_{uloc}(\Omega)$-solution of \eqref{eq:1.1} in $\Omega\times[0,T]$ if 
\begin{itemize}
  \item $u\in L^\infty(\tau,T:L^\infty(\Omega))\cap L^2(\tau,T:W^{1,2}(\Omega\cap B(0,R)))$ for any $\tau\in(0,T)$ and $R>0$,
  \item $u\in C([0,T):L^r_{uloc,\rho}(\Omega))$ with $\displaystyle{\lim_{t\to 0}}\,\|u(t)-\varphi\|_{r,\rho}=0$ for some $\rho>0$,
  \item $u$ satisfies 
  \begin{equation}
  \label{eq:1.6}
  \int_0^T\int_\Omega
  \left\{-u\partial_t\phi+\nabla u\cdot\nabla\phi\right\}\,dyds=\int_0^T\int_{\partial\Omega}|u|^{p-1}u\phi\,d\sigma ds
  \end{equation}
  for all $\phi\in C_0^\infty({\bf R}^N\times(0,T))$. 
\end{itemize}
Here $d\sigma$ is the surface measure on $\partial\Omega$. 
Furthermore, for any continuous function $u$ in $\overline{\Omega}\times(0,T)$, 
we say that $u$ is a $L^r_{uloc}(\Omega)$-solution of \eqref{eq:1.1} in $\Omega\times[0,T)$ 
if $u$ is a $L^r_{uloc}(\Omega)$-solution of \eqref{eq:1.1} in $\Omega\times[0,\eta]$ for any $\eta\in(0,T)$. 
\end{definition}
We remark the following for any $\rho$, $\rho'\in(0,\infty)$: 
\begin{itemize}
\item $f\in L^r_{uloc,\rho}(\Omega)$ is equivalent to $f\in L^r_{uloc,\rho'}(\Omega)$;
\item $u\in C([0,T]:L^r_{uloc,\rho}(\Omega))$ is equivalent to $u\in C([0,T]:L^r_{uloc,\rho'}(\Omega))$. 
\end{itemize}
These follow from property~(i) in Section~2.
\vspace{5pt}

Now we are ready to state the main  results of this paper. 
Let $p_*=1+1/N$. 
\begin{theorem}
\label{Theorem:1.1} 
Let  $N\ge 1$ and $\Omega\subset{\bf R}^N$ be a uniformly regular domain of class $C^1$.  
Let $\rho_*$ satisfy \eqref{eq:1.3} and \eqref{eq:1.5}. 
Then, for any $1\le r<\infty$ with   
\begin{equation}
\label{eq:1.7}
\left\{
\begin{array}{ll}
r\ge N(p-1) & \mbox{if}\quad p>p_*,\vspace{3pt}\\
r>1 & \mbox{if}\quad p=p_*,\vspace{3pt}\\
r\ge 1 & \mbox{if}\quad 1<p<p_*, 
\end{array}
\right.
\end{equation}
there exists a positive constant $\gamma_1$ such that, 
for any $\varphi\in {\cal L}^r_{uloc,\rho}(\Omega)$ with
\begin{equation}
\label{eq:1.8}
\rho^{\frac{1}{p-1}-\frac{N}{r}}\|\varphi\|_{r,\rho}\le\gamma_1
\end{equation}
for some $\rho\in(0,\rho_*/2)$, 
problem~\eqref{eq:1.1} possesses a $L^r_{uloc}(\Omega)$-solution $u$ of \eqref{eq:1.1} in $\Omega\times[0,\mu\rho^2]$ 
satisfying 
\begin{eqnarray}
\label{eq:1.9} 
  & &\!\!\!\!\!\!
\sup_{0<t<\mu\rho^2}\,\|u(t)\|_{r,\rho}\le C\|\varphi\|_{r,\rho},\\
  & &\!\!\!\!\!\!
 \sup_{0<t<\mu\rho^2}\,t^{\frac{N}{2r}}\|u(t)\|_{L^\infty(\Omega)}\le C\|\varphi\|_{r,\rho}. 
 \label{eq:1.10}
 \end{eqnarray}
 Here $C$ and $\mu$ are constants depending only on $N$, $\Omega$, $p$ and $r$. 
\end{theorem}
\begin{theorem}
\label{Theorem:1.2} 
Assume the same conditions as in Theorem~{\rm\ref{Theorem:1.1}}. 
Let $v$ and $w$ be $L^r_{uloc}(\Omega)$-solutions in $\Omega\times[0,T)$ 
such that $v(x,0)\le w(x,0)$ for almost all $x\in\Omega$, 
where $T>0$ and $r$ is as in \eqref{eq:1.7}. 
Assume, if $r=1$, that 
\begin{equation}
\label{eq:1.11}
\limsup_{t\to +0}\,t^{\frac{1}{2(p-1)}}\left[\|v(t)\|_{L^\infty(\Omega)}+\|w(t)\|_{L^\infty(\Omega)}\right]<\infty.
\end{equation}
Then there exists a positive constant $\gamma_2$ such that, 
if 
\begin{equation}
\label{eq:1.12}
\rho^{\frac{1}{p-1}-\frac{N}{r}}
\left[\|v(0)\|_{r,\rho}+\|w(0)\|_{r,\rho}\right]\le\gamma_2
\end{equation}
for some $\rho\in(0,\rho_*/2)$, then 
$$
v(x,t)\le w(x,t)\quad\mbox{in}\quad \Omega\times(0,T). 
$$
\end{theorem}
We give some comments related to Theorems~\ref{Theorem:1.1} and \ref{Theorem:1.2}.
\begin{itemize}
  \item[(i)] 
  Let $u$ be a $L^r_{uloc}(\Omega)$-solution of \eqref{eq:1.1} in $\Omega\times[0,T)$. 
  It follows from Definition~{\rm\ref{Definition:1.1}} that 
  $u\in L^\infty(\tau,\sigma:L^\infty(\Omega))$ for any $0<\tau<\sigma<T$. 
  This together with Theorem~{\rm 6.2} of {\rm\cite{DB01}} implies that $u(t)\in BUC(\Omega)$ 
  for any $t\in(0,T)$. This means that $u(0)\in{\cal L}^r_{uloc,\rho}(\Omega)$ for any $\rho>0$. 
  \item[(ii)] 
  Consider the case $\Omega={\bf R}^N_+$. 
  Let $u$ be a $L^r_{uloc}(\Omega)$-solution of \eqref{eq:1.1} blowing up at $t=T<\infty$, 
  where $r$ is as in \eqref{eq:1.7}.
  Then, for any $\mu>0$, 
  $u_\mu$ defined by~\eqref{eq:1.2} satisfies  
  $$
  \mu^{-\left(\frac{1}{p-1}-\frac{N}{r}\right)}\|u_\mu(0)\|_{r,\mu^{-1}}=\|u(0)\|_{r,1}
  $$
  and it blows up at $t=\mu^{-2}T$. 
  This means that  Theorem~{\rm\ref{Theorem:1.1}} holds with $\rho=1$ if and only if  
  Theorem~{\rm\ref{Theorem:1.1}} holds for any $\rho>0$. 
  \item[(iii)] 
  Let $1\le r<\infty$. If, either
  $$
  {\rm(a)}\quad f\in L^r_{uloc,1}(\Omega), \quad r>N(p-1)
  \qquad\mbox{or}\qquad
  {\rm(b)}\quad f\in L^r(\Omega),\quad r\ge N(p-1),
  $$ 
  then, for any $\gamma>0$, 
  we can find a constant $\rho>0$ such that 
  $\rho^{\frac{1}{p-1}-\frac{N}{r}}\|f\|_{r,\rho}\le\gamma$.  
\end{itemize}
As a corollary of Theorem~\ref{Theorem:1.1}, 
we have: 
\begin{corollary}
\label{Corollary:1.1}
Assume the same conditions as in Theorem~{\rm\ref{Theorem:1.1}} and $p>p_*$. 
\begin{itemize}
  \item[{\rm (i)}] 
  For any $\varphi\in L^{N(p-1)}(\Omega)$, 
  problem~\eqref{eq:1.1} has a unique $L^{N(p-1)}_{uloc}(\Omega)$-solution in $\Omega\times[0,T]$ for some $T>0$. 
  \item[{\rm (ii)}] Assume $\rho_*=\infty$. 
  Then there exists a constant $\gamma$ such that, if 
  \begin{equation}
  \label{eq:1.13}
  \|\varphi\|_{L^{N(p-1)}(\Omega)}\le\gamma,
  \end{equation}
  then problem~\eqref{eq:1.1} has a unique $L^{N(p-1)}_{uloc}(\Omega)$-solution $u$ such that 
  $$
  \sup_{0<t<\infty}\|u(t)\|_{L^{N(p-1)}(\Omega)}+\sup_{0<t<\infty}\,t^{\frac{1}{2(p-1)}}\|u(t)\|_{L^\infty(\Omega)}<\infty.
  $$
\end{itemize}
\end{corollary}
For further applications of our theorems, see Section~5. 
\begin{remark}
\label{Remark:1.1}
Let $\Omega={\bf R}^N_+:=\{(x',x_N)\in{\bf R}^N\,:\,x_N>0\}$. 
If $1<p\le p_*$, then problem~\eqref{eq:1.1} possesses no positive global-in-time solutions.   
See {\rm\cite{DFL}} and {\rm\cite{GL}}. 
For the case $p>p_*$, 
it is proved in {\rm\cite{K}} {\rm(}see also {\rm\cite{IK}}{\rm)} that,  
if $\varphi\ge 0$, $\varphi\not\equiv 0$ in $\Omega$ and 
$$
\|\varphi\|_{L^1({\bf R}^N_+)}\|\varphi\|_{L^\infty({\bf R}^N_+)}^{N(p-1)-1} 
\quad\mbox{is sufficiently small},
$$ 
then there exists a positive global-in-time solution of \eqref{eq:1.1}. 
This also immediately follows from assertion~{\rm(ii)} of Corollary~{\rm\ref{Corollary:1.1}} and the comparison principle. 
\end{remark}

\vspace{5pt}

We explain the idea of the proof of Theorem~\ref{Theorem:1.1}. 
Under the assumptions of Theorem~\ref{Theorem:1.1}, 
there exists a sequence $\{\varphi_n\}_{n=1}^\infty\subset BUC(\Omega)$ 
such that 
\begin{equation}
\label{eq:1.14}
\lim_{n\to\infty}\|\varphi-\varphi_n\|_{r,\rho}=0,
\qquad
\sup_n\,\|\varphi_n\|_{r,\rho}\le 2\|\varphi\|_{r,\rho}. 
\end{equation}
For any $n=1,2,\dots$, 
let $u_n$ satisfy in the classical sense
\begin{equation}
\label{eq:1.15}
\left\{
\begin{array}{ll}
\partial_t u=\Delta u & \mbox{in}\quad\Omega\times(0,T_n),\vspace{3pt}\\
\nabla u\cdot\nu(x)=|u|^{p-1}u & \mbox{on}\quad\partial\Omega\times(0,T_n),\vspace{3pt}\\
u(x,0)=\varphi_n(x) & \mbox{in}\quad\Omega,
\end{array}
\right.
\end{equation}
where $T_n$ is the blow-up time of the solution $u_n$. 
By regularity theorems for parabolic equations (see e.g. \cite{DB01} and \cite[Chapters~III and IV]{LSU}) 
we see that 
\begin{equation}
\label{eq:1.16}
u_n\in BUC(\overline{\Omega}\times[0,T]),
\qquad
\nabla u_n\in L^\infty(\Omega\times(\tau,T)),
\end{equation}
for any $0<\tau<T<T_n$, which imply that 
$u_n$ is a $L^r_{uloc}(\Omega)$-solution in $\Omega\times[0,T_n)$ for any $1\le r<\infty$.  
Set 
$$
\Psi_{r,\rho}[u_n](t):=\sup_{0\le\tau\le t}\,\sup_{x\in\overline{\Omega}}\,
\int_{\Omega(x,\rho)}|u_n(y,\tau)|^r\,dy,
\qquad 0\le t<T_n.
$$
It follows from \eqref{eq:1.8} and \eqref{eq:1.14} that 
\begin{equation}
\label{eq:1.17}
\Psi_{r,\rho}[u_n](0)^{\frac{1}{r}}=\|\varphi_n\|_{r,\rho}\le 2\|\varphi\|_{r,\rho}
\le 2\gamma_1\rho^{-\frac{1}{p-1}+\frac{N}{r}}.
\end{equation}
Define 
\begin{equation}
\label{eq:1.18}
\begin{split}
T_n^* & :=\sup\left\{\sigma\in(0,T_n)\,:\,\Psi_{r,\rho}[u_n](t)\le 6M\Psi_{r,\rho}[u_n](0)\quad\mbox{in}\quad[0,\sigma]\right\},\\
T_n^{**} & :=\sup\left\{\sigma\in(0,T_n)\,:\,
\rho^{-1}+\|u_n(t)\|_{L^\infty(\Omega)}^{p-1}\le 2t^{-\frac{1}{2}}\quad\mbox{in}\quad(0,\sigma]\right\},
\end{split}
\end{equation}
where $M$ is the integer given in Lemma~\ref{Lemma:2.1}. 
We adapt the arguments in \cite{A}, \cite{AD} and \cite{I} to 
obtain uniform estimates of $u_n$ and $u_m-u_n$ with respect to $m$, $n=1,2,\dots$, 
and prove that 
$$
\inf_nT_n^*\ge\mu\rho^2,\qquad \inf_nT_n^{**}\ge\mu\rho^2,
$$ 
for some $\mu>0$. This enables us to prove Theorem~\ref{Theorem:1.1}. 
Theorem~\ref{Theorem:1.2} follows from a similar argument as in Theorem~\ref{Theorem:1.1}. 
\vspace{5pt}

The rest of this paper is organized as follows. 
In Section 2 we give some preliminary lemmas related to $\rho_*$. 
In Sections~3 and 4 we prove Theorems~\ref{Theorem:1.1} and \ref{Theorem:1.2}. 
In Section 5, as applications of Theorem~\ref{Theorem:1.1}, 
we give some results on the blow-up time and the blow-up rate of the solutions. 
\section{Preliminaries}
In this section we recall some properties of uniformly local $L^r$ spaces 
and prove some lemmas related to $\rho_*$. 
Furthermore, we give some inequalities used in Sections 3 and 4. 
In what follows, 
the letter $C$ denotes a generic constant independent of $x\in\overline{\Omega}$, $n$ and $\rho$. 
\vspace{5pt}

Let $1\le r<\infty$. We first recall the following properties of $L^r_{uloc,\rho}(\Omega)$: 
\begin{itemize}
  \item[(i)] 
  if $f\in L^r_{uloc,\rho}(\Omega)$ for some $\rho>0$, 
  then, for any $\rho'>0$,  $f\in L^r_{uloc,\rho'}(\Omega)$ and 
  $$
  \|f\|_{r,\rho'}\le C_1\|f\|_{r,\rho}
  $$ 
  for some constant $C_1$ depending only on $N$, $\rho$ and $\rho'$;
  \item[(ii)] there exists a constant $C_2$ depending only on $N$ such that
  \begin{equation}
  \label{eq:2.1}
  \|f\|_{r,\rho}\le C_2\rho^{N(\frac{1}{r}-\frac{1}{q})}\|f\|_{q,\rho},
  \qquad
  f\in L^q_{uloc,\rho}(\Omega), 
  \end{equation}
  for any $1\le r\le q<\infty$ and $\rho>0$;
  \item[(iii)] if $f\in L^r(\Omega)$, then 
  $f\in L^r_{uloc,\rho}(\Omega)$ for any $\rho>0$ and 
  \begin{equation}
  \label{eq:2.2}
  \lim_{\rho\to +0}\,\|f\|_{r,\rho}=0. 
  \end{equation}
\end{itemize}
Properties~(ii) and (iii) are proved by the H\"older inequality 
and the absolute continuity of $|f|^r\,dy$ with respect to $dy$. 
Property~(i) follows from the following lemma. 
\begin{lemma}
\label{Lemma:2.1}
Let $N\ge 1$ and $\Omega$ be a domain in ${\bf R}^N$. 
Then there exists $M\in\{1,2,\dots\}$ depending only on $N$ such that, 
for any $x\in\overline{\Omega}$ and $\rho>0$, 
\begin{equation}
\label{eq:2.3}
\Omega(x,2\rho)\subset\bigcup_{k=1}^n \Omega(x_k,\rho)
\end{equation}
for some $\{x_k\}_{k=1}^n\subset\overline{\Omega}$ with $n\le M$. 
\end{lemma}
{\bf Proof.}
There exist $M\in\{1,2,\dots\}$ and $\{y_k\}_{k=1}^M\subset B(0,2)$ such that  
$$
B(0,2)\subset\bigcup_{k=1}^M B(y_k,1/2). 
$$
Then, for any $x\in\overline{\Omega}$ and $\rho>0$, 
we can find $\{y_{k_i}\}_{i=1}^n\subset\{y_k\}_{k=1}^M$ such that 
\begin{equation}
\label{eq:2.4}
\Omega(x+\rho y_{k_i},\rho/2)\not=\emptyset
\quad\mbox{and}\quad
\Omega(x,2\rho)\subset\bigcup_{i=1}^n \Omega(x+\rho y_{k_i},\rho/2). 
\end{equation}
Furthermore, for any $i\in\{1,\dots,n\}$, 
there exists $x_{k_i}\in\overline{\Omega}$ such that 
$$
x_{k_i}\in \Omega(x+\rho y_{k_i},\rho/2)
\quad\mbox{and}
\quad
\Omega(x+\rho y_{k_i},\rho/2)\subset \Omega(x_{k_i},\rho).
$$ 
This together with \eqref{eq:2.4} implies \eqref{eq:2.3}, and Lemma~\ref{Lemma:2.1} follows.
$\Box$
\vspace{7pt}

We state a lemma on the existence of $\rho_*$ satisfying \eqref{eq:1.3} and \eqref{eq:1.5}. 
\begin{lemma}
\label{Lemma:2.2}
Let $N\ge 1$ and $\Omega$ be a uniformly regular domain of class $C^1$.  
Then there exists $\rho_*>0$ such that \eqref{eq:1.3} and \eqref{eq:1.5} hold. 
In particular, if 
\begin{equation}
\label{eq:2.5}
\Omega=\{(x',x_N)\in{\bf R}^N\,:\,x_N>\Phi(x')\},
\end{equation}
where $N\ge 2$ and $\Phi\in C^1({\bf R}^{N-1})$ with $\|\nabla\Phi\|_{L^\infty({\bf R}^{N-1})}<\infty$, 
then \eqref{eq:1.3} and \eqref{eq:1.5} hold with $\rho_*=\infty$. 
\end{lemma}
{\bf Proof.}
By the definition of uniformly regular domain, it suffices to consider the case \eqref{eq:2.5}. 
Let $f\in C^1_0(B(x_*,\rho))$, where $x_*\in\overline{\Omega}$ and $\rho>0$. 
Set $f=0$ outside $B(x_*,\rho)$.  
We first consider the case of $\partial\Omega(x_*,\rho)\not=\emptyset$. 
Then there exists $y_*\in\partial\Omega$ such that $B(x_*,\rho)\subset B(y_*,2\rho)$. 
Set 
$$
g(x',x_N):=
\left\{
\begin{array}{ll}
f(x'-y_*',x_N+\Phi(x')) & \mbox{for}\quad x_N\ge 0,\vspace{3pt}\\
f(x'-y_*',-x_N+\Phi(x')) & \mbox{for}\quad x_N<0,
\end{array}
\right.
\qquad
\tilde{g}(z):=g(2\rho'z),
$$
where 
$$
\rho'=\rho\left(1+\|\nabla\Phi\|_{L^\infty({\bf R}^{N-1})}^2\right)^{1/2}.
$$
Then $\tilde{g}\in C^1_0(B(0,1))$. 
Applying the Gagliardo-Nirenberg  inequality (see e.g. \cite{GGS})
and the trace imbedding theorem (see e.g. \cite[Theorem~5.22]{A}), 
we obtain 
\begin{eqnarray*}
 \|\tilde{g}\|_{L^\beta(B(0,1))} \!\!\! & \le &\!\!\! C\|\tilde{g}\|_{L^\beta(B(0,1))}^{1-\sigma}\|\nabla \tilde{g}\|_{L^2(B(0,1))}^\sigma,\\
\int_{B(0,1)\cap\partial{\bf R}^N_+}|\tilde{g}|\,d\sigma
 \!\!\! & \le &\!\!\!  
 C\|\tilde{g}\|_{W^{1,1}(B(0,1)\cap {\bf R}^N_+)}
 \le C\|\nabla\tilde{g}\|_{L^1(B(0,1)\cap{\bf R}^N_+)},
\end{eqnarray*}
where $\alpha$, $\beta$ and $\sigma$ are as in \eqref{eq:1.4} and $\alpha\not=\infty$ if $N=2$. 
These imply that 
\begin{eqnarray*}
 \|g\|_{L^\beta(B(0,2\rho'))}
 \!\!\! & \le &\!\!\! 
 C\|g\|_{L^\beta(B(0,2\rho'))}^{1-\sigma}\|\nabla g\|_{L^2(B(0,2\rho'))}^\sigma,\\
\int_{B(0,2\rho')\cap\partial{\bf R}^N_+}|g|\,d\sigma
 \!\!\! & \le &\!\!\!  C\|\nabla g\|_{L^1(B(0,2\rho')\cap{\bf R}^N_+)},
\end{eqnarray*}
for some constants $C$ independent of $\rho$. 
Then we have 
\begin{eqnarray}
\|f\|_{L^\beta(\Omega(x_*,\rho))} & = & \|f\|_{L^\beta(\Omega(y_*,2\rho))}
\le C\|g\|_{L^\beta(B(0,2\rho'))}\\
 & \le & C\|g\|_{L^\beta(B(0,2\rho'))}^{1-\sigma}\|\nabla g\|_{L^2(B(0,2\rho'))}^\sigma\notag\\
 & \le & C\|f\|_{L^\beta(\Omega(x_*,\rho))}^{1-\sigma}\|\nabla f\|_{L^2(\Omega(x_*,\rho))}^\sigma,\label{eq:2.7}\\
\int_{\partial\Omega(x_*,\rho)}|f|\,d\sigma
 & \le & C\int_{B(0,2\rho')\cap\partial{\bf R}^N_+}|g|\,d\sigma
 \le C\|\nabla g\|_{L^1(B(0,2\rho')\cap{\bf R}^N_+)}\notag\\
 & \le & C\|\nabla f\|_{L^1(\Omega(x_*,\rho))}. 
\label{eq:2.8}
\end{eqnarray}
Therefore we obtain \eqref{eq:1.3} and \eqref{eq:1.5} for any $\rho>0$ 
in the case of $\partial\Omega(x_*,\rho)\not=\emptyset$. 
Similarly, we get \eqref{eq:1.3} and \eqref{eq:1.5} for all $\rho>0$ 
in the case of $\partial\Omega(x_*,\rho)=\emptyset$. 
Thus \eqref{eq:1.3} and \eqref{eq:1.5} hold with $\rho_*=\infty$ in the case \eqref{eq:2.5}, 
and the proof is complete.  
$\Box$\vspace{5pt}

We obtain the following two lemmas by using \eqref{eq:1.3} and \eqref{eq:1.5}. 
\begin{lemma}
\label{Lemma:2.3}
Let $N\ge 1$ and $\Omega\subset{\bf R}^N$ be a uniformly regular domain of class $C^1$.  
Let $\rho_*$ satisfy \eqref{eq:1.3} and \eqref{eq:1.5}. 
Then there exists a constant $C_1$ such that 
\begin{equation}
\label{eq:2.9}
\int_{\partial\Omega(x,\rho)}\phi^2\,d\sigma
\le\epsilon\int_{\Omega(x,\rho)}|\nabla\phi|^2\,dy+\frac{C_1}{\epsilon}\int_{\Omega(x,\rho)}\phi^2\,dy
\end{equation}
for all $\phi\in C^1_0(B(x,\rho))$, $\epsilon>0$, $x\in\overline{\Omega}$ and $\rho\in(0,\rho_*)$. 
Furthermore, for any $p>1$ and $r>0$, 
there exists a constant $C_2$ such that 
\begin{equation}
\label{eq:2.10}
\int_{\Omega(x,\rho)}\,f^{2p+r-2}\,dy
\le C_2\left(\int_{\Omega(x,\rho)}\,f^{N(p-1)}\,dy\right)^{\frac{2}{N}}
\int_{\Omega(x,\rho)}\,|\nabla f^{\frac{r}{2}}|^2\,dy
\end{equation}
for all nonnegative functions $f$ satisfying $f^{r/2}\in C^1(\Omega(x,\rho))$ with $f=0$ near $\Omega\,\cap\,\partial B(x,\rho)$, 
$\rho\in(0,\rho_*)$ and $x\in\overline{\Omega}$. 
\end{lemma}
{\bf Proof.}
It follows from \eqref{eq:1.5} that 
\begin{equation*}
\begin{split}
\int_{\partial\Omega(x,\rho)}\phi^2\,d\sigma
 & \le C\int_{\Omega(x,\rho)}|\nabla \phi^2|\,dy
\le 2C\int_{\Omega(x,\rho)}|\phi|\,|\nabla \phi|\,dy\\
 & \le\epsilon\int_{\Omega(x,\rho)}|\nabla\phi|^2\,dy
 +\frac{C^2}{\epsilon}\int_{\Omega(x,\rho)}\phi^2\,dy
\end{split}
\end{equation*}
for all $\phi\in W^{1,2}_0(B(x,\rho))$, $\epsilon>0$, $x\in\overline{\Omega}$ and $\rho\in(0,\rho_*)$. 
This implies \eqref{eq:2.9}. 

Let $r>0$ and $0<\rho<\rho_*$. 
If $2N(p-1)\ge r$, 
then, by \eqref{eq:1.5} we have 
\begin{equation}
\label{eq:2.11}
\int_{\Omega(x,\rho)}\,g^{\frac{4}{r}(p-1)+2}\,dy
\le C\left(\int_{\Omega(x,\rho)}\,g^{\frac{2N(p-1)}{r}}\,dy\right)^{\frac{2}{N}}
\int_{\Omega(x,\rho)}\,|\nabla g|^2\,dy
\end{equation}
for all $g\in C^1_0(B(x,\rho))$ and $x\in\overline{\Omega}$. 
Furthermore, we obtain \eqref{eq:2.11} by the H\"older inequality and \eqref{eq:1.5} even for the case $2N(p-1)<r$ 
(see e.g. \cite[Lemma~3]{N}).   
Then, setting $g=f^{r/2}$, we obtain \eqref{eq:2.10}, and the proof is complete. 
$\Box$ 
\begin{lemma}
\label{Lemma:2.4}
Assume the same conditions as in Theorem~{\rm\ref{Theorem:1.1}}. 
Let $r\ge 1$, $T>0$ and $f$ be a nonnegative function such that 
$$
f\in C([0,T]:L^r_{uloc,\rho}(\Omega))\cap L^2(\tau,T:W^{1,2}(\Omega\cap B(0,R)))
$$ 
for any $\rho\in(0,\rho_*/2)$, $\tau\in(0,T)$ and $R>0$. 
Let $x\in\overline{\Omega}$ and $\zeta$ be a smooth function in ${\bf R}^N$ such that 
\begin{eqnarray*}
 & & 0\le\zeta\le 1\quad\mbox{and}\quad|\nabla\zeta|\le 2\rho^{-1}\quad\mbox{in}\quad{\bf R}^N,\vspace{5pt}\\
 & & \zeta=1\quad\mbox{on}\quad B(x,\rho),
\quad
\zeta=0\quad\mbox{outside}\quad B(x,2\rho).
\end{eqnarray*} 
Set $f_\epsilon=f+\epsilon$ for $\epsilon>0$. 
Then, for any sufficiently large $k\ge 2$, 
there exists a constant $C$ such that 
\begin{equation}
\label{eq:2.12}
\begin{split}
 & \sup_{x\in\overline{\Omega}}\,\int_\tau^t\int_{\partial\Omega(x,2\rho)}f_\epsilon^{p+r-1}\zeta^k\,d\sigma ds\\
 & \le C\biggr[\rho^{\frac{r}{p-1}-N}\Psi_{r,\rho}[f_\epsilon](t)\biggr]^{\frac{p-1}{r}}
\left[\sup_{x\in\overline{\Omega}}\,\int_\tau^t\int_{\Omega(x,\rho)}|\nabla f_\epsilon^{\frac{r}{2}}|^2\,dyds
+\rho^{-2}(t-\tau)\Psi_{r,\rho}[f_\epsilon](t)\right]
\end{split}
\end{equation}
for all $0<\tau<t\le T$, $\rho\in(0,\rho_*/2)$ and $\epsilon>0$. 
\end{lemma}
{\bf Proof.}
Let $\rho\in(0,\rho_*/2)$. 
It suffices to consider the case where $\partial\Omega(x,\rho)\not=\emptyset$. 
Let $k\ge 2$ be such that 
\begin{equation}
\label{eq:2.13}
\frac{k}{2p+r-2}\cdot\frac{r}{2}\ge 1. 
\end{equation}
By \eqref{eq:1.3} and Lemma~\ref{Lemma:2.1}, 
for any $\delta>0$, we have  
\begin{equation}
\label{eq:2.14}
\begin{split}
 & \int_\tau^t\int_{\partial\Omega(x,2\rho)}f_\epsilon^{p+r-1}\zeta^k\,d\sigma ds
 \le C\int_\tau^t\int_{\Omega(x,2\rho)}\left|\nabla[f_\epsilon^{p+r-1}\zeta^k]\right|\,dyds\\
 & \le C\int_\tau^t\int_{\Omega(x,2\rho)}f_\epsilon^{p+\frac{r}{2}-1}|\nabla f_\epsilon^{\frac{r}{2}}|\zeta^k\,dyds
 +C\int_\tau^t\int_{\Omega(x,2\rho)}f_\epsilon^{p+r-1}|\nabla\zeta|\zeta^{k-1}\,dyds\\
 & \le C\delta\int_\tau^t\int_{\Omega(x,2\rho)}f_\epsilon^{2p+r-2}\zeta^k\,dyds\\
 & \qquad\quad
+C\delta^{-1}\int_\tau^t\int_{\Omega(x,2\rho)}|\nabla f_\epsilon^{\frac{r}{2}}|^2\zeta^k\,dyds
+C\delta^{-1}\int_\tau^t\int_{\Omega(x,2\rho)}f_\epsilon^r\zeta^{k-2}|\nabla\zeta|^2\,dy\,ds\\
 & \le C\delta\int_\tau^t\int_{\Omega(x,2\rho)}f_\epsilon^{2p+r-2}\zeta^k\,dyds\\
 & \qquad\quad
+C\delta^{-1}\sup_{x\in\overline{\Omega}}\,\int_\tau^t\int_{\Omega(x,\rho)}|\nabla f_\epsilon^{\frac{r}{2}}|^2\,dyds
+C\delta^{-1}\rho^{-2}(t-\tau)\Psi_{r,\rho}[f_\epsilon](t)
\end{split}
\end{equation}
for $0<\tau<t\le T$, 
where $C$ is a constant independent of $\epsilon$ and $\delta$. 
Set $g_\epsilon:=f_\epsilon\zeta^{k/(2p+r-2)}$. 
It follows from \eqref{eq:2.13} that 
$f_\epsilon^{r/2}=0$ near $\partial B(x,2\rho)\cap\Omega$. 
Then, by Lemmas~\ref{Lemma:2.1} and \ref{Lemma:2.3} 
we have 
\begin{equation}
\label{eq:2.15}
\begin{split}
 & \int_\tau^t\int_{\Omega(x,2\rho)}f_\epsilon(y,\tau)^{2p+r-2}\zeta^k\,dyds
 =\int_\tau^t\int_{\Omega(x,2\rho)}g_\epsilon(y,\tau)^{2p+r-2}\,dyds\\
 & \le C\sup_{0<s<t}\left(\int_{\Omega(x,2\rho)}g_\epsilon(y,s)^{N(p-1)}\,dy\right)^{\frac{2}{N}}
\int_\tau^t\int_{\Omega(x,2\rho)}|\nabla g_\epsilon^{\frac{r}{2}}|^2\,dyds\\
 & \le C\sup_{0<s<t}\left(\rho^{\frac{r}{p-1}-N}\int_{\Omega(x,2\rho)}f_\epsilon(y,s)^r\,dy\right)^{\frac{2(p-1)}{r}}\\
 & \qquad\qquad\qquad
\times\biggr[\int_\tau^t\int_{\Omega(x,2\rho)}|\nabla f_\epsilon^{\frac{r}{2}}|^2\,dyds
+\rho^{-2}\int_\tau^t\int_{\Omega(x,2\rho)}f_\epsilon^r\,dyds\biggr]\\
 & \le C\biggr[\rho^{\frac{r}{p-1}-N}\Psi_{r,\rho}[f_\epsilon](t)\biggr]^{\frac{2(p-1)}{r}}\\
 & \qquad\qquad\qquad
\times\biggr[\sup_{x\in\overline{\Omega}}\int_\tau^t\int_{\Omega(x,\rho)}|\nabla f_\epsilon^{\frac{r}{2}}|^2\,dyds
+\rho^{-2}(t-\tau)\Psi_{r,\rho}[f_\epsilon](t)\biggr]
\end{split}
\end{equation}
for $0<\tau<t\le T$. 
Therefore, taking 
$\delta=[\rho^{\frac{r}{p-1}-N}\Psi_{r,\rho}[f_\epsilon](t)]^{-(p-1)/r}$,  
by \eqref{eq:2.14} and \eqref{eq:2.15} we obtain \eqref{eq:2.12}, 
and the proof is complete. 
$\Box$
\section{Proof of Theorems~\ref{Theorem:1.1} and \ref{Theorem:1.2} in the case $r>1$.}
Let $v$ and $w$ be $L^r_{uloc}(\Omega)$-solutions 
of \eqref{eq:1.1} in $\Omega\times[0,T]$, where $0<T<\infty$ and $r$ is as in \eqref{eq:1.7}. 
Set $z:=v-w$ and $z_\epsilon:=\max\{z,0\}+\epsilon$ for $\epsilon\ge 0$. 
Then $z_\epsilon$ satisfies 
\begin{equation}
\label{eq:3.1}
\partial_t z_\epsilon\le\Delta z_\epsilon\quad\mbox{in}\quad\Omega\times(0,T],\qquad
\nabla z_\epsilon\cdot\nu(x)\le a(x,t)z_\epsilon\quad\mbox{on}\quad\partial\Omega\times(0,T],
\end{equation}
in the weak sense (see e.g. \cite[Chapter~II]{DB02}). 
Here 
\begin{equation}
\label{eq:3.2}
a(x,t):=\left\{
\begin{array}{ll}
\displaystyle{\frac{|v(x,t)|^{p-1}v(x,t)-|w(x,t)|^{p-1}w(x,t)}{v(x,t)-w(x,t)}} & \mbox{if}\quad v(x,t)\not=w(x,t),\vspace{5pt}\\
p|v(x,t)|^{p-1} & \mbox{if}\quad v(x,t)=w(x,t), 
\end{array}
\right.
\end{equation}
which satisfies 
\begin{equation}
\label{eq:3.3}
0\le a(x,t)\le C(|v|^{p-1}+|w|^{p-1})
\quad\mbox{in}\quad\Omega\times(0,T]. 
\end{equation}
In this section we give some estimates of $z$, 
and prove Theorems~\ref{Theorem:1.1} and \ref{Theorem:1.2} in the case $r>1$. 

We first give an $L^\infty_{loc}$ estimate of $z_0$
by using the Moser iteration method with the aid of \eqref{eq:1.18}. 
For related results, see \cite{FiloK}. 
\begin{lemma}
\label{Lemma:3.1}
Assume the same conditions as in Theorem~{\rm\ref{Theorem:1.1}}. 
Let $v$ and $w$ be $L^r_{uloc}(\Omega)$-solutions 
of \eqref{eq:1.1} in $\Omega\times[0,T]$, where $0<T<\infty$ and $r\ge 1$. 
Set $z_0:=\max\{v-w,0\}$ and $a=a(x,t)$ as in \eqref{eq:3.2}. 
Then there exists a constant $C$ such that 
\begin{eqnarray}
\label{eq:3.4}
 & & 
\|z_0(t)\|_{L^\infty(\Omega(x,R_1)\times(t_1,t))}
\le CD^{\frac{N+2}{2r}}\left(\int_{t_2}^t\int_{\Omega(x,R_2)}\,z_0^r\,dyds\right)^{1/r},\\
\label{eq:3.5}
 & & \int_{t_1}^t\int_{\Omega(x,R_1)}|\nabla z_0|^2\,dyds
\le CD\int_{t_2}^t\int_{\Omega(x,R_2)} z_0^2\,dyds, 
\end{eqnarray}
for all $x\in\overline{\Omega}$, $0<R_1<R_2<\rho_*$ and $0<t_2<t_1<t\le T$, 
where 
$$
D:=\|a\|_{L^\infty(\Omega(x,R_2)\times(t_2,t))}^2+(R_2-R_1)^{-2}+(t_1-t_2)^{-1}.
$$
\end{lemma}
{\bf Proof.}
Let $x\in\overline{\Omega}$, $0<R_1<R_2<\rho_*$ and $0<t_2<t_1<t\le T$. 
For $j=0,1,2,\dots$, 
set
$$
r_j:=R_1+(R_2-R_1)2^{-j},\quad
\tau_j:=t_1-(t_1-t_2)2^{-j},\quad
Q_j:=\Omega(x,r_j)\times(\tau_j,t). 
$$
Let $\zeta_j$ be a piecewise smooth function in $Q_j$ such that
\begin{equation}
\begin{split}
 & 0\le\zeta_j\le 1\quad \mbox{in} \quad {\bf R}^N,
\quad \zeta_j=1\quad\mbox{on}\quad Q_{j+1},\vspace{3pt}\\
 & \zeta_j=0\quad\mbox{near}\quad \partial\Omega(x,r_j)\times[\tau_j,t]\cup\Omega(x,r_j)\times\{\tau_j\},\\
 & |\nabla\zeta_j|\le\frac{2^{j+1}}{R_2-R_1}
\quad\mbox{and}\quad
0\le\partial_t\zeta_j\le\frac{2^{j+1}}{t_1-t_2}
\quad\mbox{in}\quad Q_j.
\end{split}
\label{eq:3.6}
\end{equation}
Let $\alpha_0>1$ and $\epsilon>0$. 
For any $\alpha\ge\alpha_0$, 
multiplying \eqref{eq:3.1} by $z_\epsilon^{\alpha-1}\zeta_j^2$ and integrating it on $Q_j$, 
we obtain 
\begin{equation}
\begin{split}
 & \frac{1}{\alpha}\sup_{\tau_j<s<t}\int_{\Omega(x,r_j)}z_\epsilon^\alpha\zeta_j^2 \,dy
+\frac{\alpha-1}{2}\iint_{Q_j}z_\epsilon^{\alpha-2}|\nabla z_\epsilon|^2\zeta_j^2\,dyds\\
 & 
\le\frac{4}{\alpha}\iint_{Q_j}z_\epsilon^\alpha\zeta_j|\partial_t\zeta_j|\,dyds
+\frac{4}{\alpha-1}\iint_{Q_j}z_\epsilon^\alpha|\nabla\zeta_j|^2\,dyds\\
 & \qquad\qquad\qquad\qquad\qquad\qquad\quad
 +2\int_{\tau_j}^t\int_{\partial\Omega(x,r_j)}a(y,s)z_\epsilon^\alpha \zeta_j^2\,d\sigma ds. 
\end{split}
\label{eq:3.7}
\end{equation}
This calculation is somewhat formal, however it is justified 
by the same argument as in~\cite[Chapter~III]{LSU} (see also \cite{DB02}). 
Then it follows that 
\begin{equation}
\begin{split}
 & \sup_{\tau_j<s<t}\int_{\Omega(x,r_j)}z_\epsilon^\alpha\zeta_j^2 \,dy
+\iint_{Q_j}|\nabla[z_\epsilon^{\frac{\alpha}{2}}\zeta_j]|^2\,dyds
\le C\iint_{Q_j}z_\epsilon^\alpha\zeta_j\partial_t\zeta_j\,dyds\\
 & \qquad\qquad\qquad
+C\iint_{Q_j}z_\epsilon^\alpha|\nabla \zeta_j|^2 \,dyds
+C\alpha\int_{\tau_j}^t\int_{\partial\Omega(x,r_j)}a(y,s)z_\epsilon^\alpha\zeta_j^2\,d\sigma ds
\end{split}
\label{eq:3.8}
\end{equation}
for all $j=0,1,2,\dots$ and $\alpha\ge\alpha_0$.  
On the other hand, 
by Lemma~\ref{Lemma:2.3} 
we have 
\begin{equation}
\begin{split}
 & C\alpha\int_{\tau_j}^t\int_{\partial\Omega(x,r_j)}a(y,s)z_\epsilon^\alpha\zeta_j^2\,d\sigma ds
\le C\alpha\|a\|_{L^\infty(Q_0)}\int_{\tau_j}^t\int_{\partial\Omega_j}
z_\epsilon^\alpha\zeta_j^2\,d\sigma ds\\
 & \qquad\quad
\le\frac{1}{2}\iint_{Q_j}
|\nabla[z_\epsilon^{\frac{\alpha}{2}}\zeta_j]|^2\,dyds
+C\alpha^2\|a\|_{L^\infty(Q_0)}^2\iint_{Q_j}z_\epsilon^\alpha\zeta_j^2\,dyds.
\end{split}
\label{eq:3.9}
\end{equation}
We deduce from \eqref{eq:3.6}, \eqref{eq:3.8} and \eqref{eq:3.9} that
\begin{equation}
\begin{split}
& \sup_{\tau_j<s<t}\int_{\Omega(x,r_j)}z_\epsilon^\alpha\zeta_j^2 \,dy
+\iint_{Q_j}|\nabla[z_\epsilon^{\frac{\alpha}{2}}\zeta_j]|^2\,dyds\\
 & \le C\left[\alpha^2\|a\|_{L^\infty(Q_0)}^2+\frac{2^{2j}}{(R_2-R_1)^2}+\frac{2^j}{t_1-t_2}\right]
\iint_{Q_j}z_\epsilon^\alpha\,dyds
\end{split}
\label{eq:3.10}
\end{equation}
for all $j=0,1,2,\dots$ and $\alpha\ge \alpha_0$. 
This together with \eqref{eq:1.5} implies that 
\begin{equation}
\begin{split}
 & \left(\iint_{Q_{j+1}}z_\epsilon^{\kappa\alpha}\,dyds\right)^{1/\kappa}\\
 & \le C\left[\alpha^2\|a\|_{L^\infty(Q_0)}^2+\frac{2^{2j}}{(R_2-R_1)^2}+\frac{2^j}{t_1-t_2}\right]
\iint_{Q_j}z_\epsilon^\alpha\,dyds
\end{split}
\label{eq:3.11}
\end{equation}
for all $j=0,1,2,\dots$ and $\alpha\ge\alpha_0$, where $\kappa:=1+2/N$. 
Furthermore, by \eqref{eq:3.10} with $\alpha=2$ 
we have \eqref{eq:3.5}. 

We prove \eqref{eq:3.4} in the case $r\ge 2$. Setting
$$
I_j:=\|z_\epsilon\|_{L^{\alpha_j}(Q_j)},
\qquad
\alpha_j:=r\kappa^j,
$$
by \eqref{eq:3.11} we have 
\begin{equation}
\label{eq:3.12}
I_{j+1}\le C^{\frac{1}{\alpha_j}}
\left[\alpha_j^2\|a\|_{L^\infty(Q_0)}^2+\frac{2^{2j}}{(R_2-R_1)^2}+\frac{2^j}{t_1-t_2}\right]^{\frac{1}{\alpha_j}} I_j
\le C^{\frac{j}{\alpha_j}}(CD)^{\frac{1}{\alpha_j}}I_j
\end{equation}
for all $j=0,1,2,\dots$, where $D:=\|a\|_{L^\infty(Q_0)}^2+(R_2-R_1)^{-2}+(t_1-t_2)^{-1}$. 
Since
$$
\sum_{j=0}^\infty\frac{1}{\alpha_j}=\frac{1}{r}\sum_{j=0}^\infty\kappa^{-j}
=\frac{1}{r(1-\kappa^{-1})}=\frac{N+2}{2r},
\qquad
\sum_{j=0}^\infty \frac{j}{\alpha_j}<\infty,
$$
we deduce from \eqref{eq:3.12} that 
$$
\|z_\epsilon\|_{L^\infty(Q_\infty)}=
\lim_{j\to\infty}I_j\le C^{\sum_{j=0}^\infty\frac{j}{\alpha_j}}(CD)^{\sum_{j=0}^\infty\frac{1}{\alpha_j}}I_0\\
\le CD^{(N+2)/2r}\|z_\epsilon\|_{L^r(Q_0)},
$$
which implies 
\begin{equation}
\label{eq:3.13}
\|z_\epsilon\|_{L^\infty(\Omega(x,R_1)\times(t_1,t))}\le
CD^{\frac{N+2}{2r}}\left(\int_{t_2}^t\int_{\Omega(x,R_2)}z_\epsilon^r \,dyds\right)^{1/r}, 
\end{equation}
where $r\ge 2$. 
Then, passing the limit as $\epsilon\to 0$, we obtain \eqref{eq:3.5}. 

On the other hand, for the case $1\le r<2$, 
applying \eqref{eq:3.13} with $r=2$ to the cylinders $Q_j$ and $Q_{j+1}$, 
we have 
\begin{equation*}
\begin{split}
\|z_\epsilon\|_{L^\infty(Q_{j+1})} & \le C\left((2^{2j}D)^{\frac{N+2}{2}}\iint_{Q_j}z_\epsilon^2\,dyds\right)^{\frac{1}{2}}\\
 & \le Cb^j\|z_\epsilon\|_{L^\infty(Q_j)}^{1-r/2}\left(D^{(N+2)/2}\iint_{Q_j}z_\epsilon^r\,dyds\right)^{\frac{1}{2}},
\end{split}
\end{equation*}
where $b=2^{(N+2)/2}$. 
Then, for any $\nu>0$, we have 
\begin{equation*}
\begin{split}
\|z_\epsilon\|_{L^\infty(Q_{j+1})}
 & \le\nu\|z_\epsilon\|_{L^\infty(Q_j)}+C\nu^{-\frac{2-r}{r}}
 b^{\frac{2}{r}j}D^{\frac{N+2}{2r}}\biggr(\iint_{Q_j}z_\epsilon^r\,dyds\biggr)^{1/r}\\
 & \le\nu^{j+1}\|z_\epsilon\|_{L^\infty(Q_0)}
 +C\nu^{-\frac{2-r}{r}}\sum_{i=0}^j(\nu b^{\frac{2}{r}})^iD^{\frac{N+2}{2r}}\biggr(\iint_{Q_0}z_\epsilon^r\,dyds\biggr)^{1/r}
\end{split}
\end{equation*}
for $j=1,2,\dots$. 
Taking a sufficiently small $\nu$ if necessary, 
we see that 
$$
\|z_\epsilon\|_{L^\infty(Q_{j+1})}\le\nu^{j+1}\|z_\epsilon\|_{L^\infty(Q_0)}+CD^{\frac{N+2}{2r}}\biggr(\iint_{Q_0}z_\epsilon^r\,dyds\biggr)^{1/r}
$$
for $j=1,2,\dots$. 
Passing to the limit as $j\to\infty$ and $\epsilon\to 0$, 
we obtain 
$$
\|z_0\|_{L^\infty(Q_\infty)}\le CD^{\frac{N+2}{2r}}\biggr(\iint_{Q_0}z_0^r\,dyds\biggr)^{1/r},
$$
which implies \eqref{eq:3.5} in the case $1\le r<2$. 
Thus Lemma~\ref{Lemma:3.1} follows.  
$\Box$
\begin{lemma}
\label{Lemma:3.2}
Assume the same conditions as in Theorem~{\rm\ref{Theorem:1.1}}. 
Let $r$ satisfy \eqref{eq:1.7} and $r>1$. 
Let $v$ be a $L^r_{uloc}(\Omega)$-solution of \eqref{eq:1.1} in $\Omega\times[0,T]$, where $T>0$. 
Then there exists a positive constant $\Lambda$ such that, 
if 
\begin{equation}
\label{eq:3.14}
\rho^{\frac{r}{p-1}-N}\Psi_{r,\rho}[v](T)\le\Lambda
\end{equation}
for some $\rho\in(0,\rho_*/2)$, then 
\begin{eqnarray}
\label{eq:3.15}
 & & \Psi_{r,\rho}[v](t)\le 5M\Psi_{r,\rho}[v](\tau),\\
\label{eq:3.16}
 & & \sup_{x\in\overline{\Omega}}\,\int_\tau^t\int_{\partial\Omega(x,\rho)}|v|^{p+r-1}\,d\sigma ds
 \le C\Lambda^{\frac{p-1}{r}}\Psi_{r,\rho}[v](\tau),
\end{eqnarray}
for all $0\le \tau\le t\le T$ with $t-\tau\le\mu\rho^2$, 
where $C$ and $\mu$ are positive constants depending only on $N$, $\Omega$, $p$ and $r$. 
\end{lemma}
{\bf Proof.} 
Let $x\in\overline{\Omega}$ and let $\zeta$ and $k$ be as in Lemma~\ref{Lemma:2.4}. 
By \eqref{eq:3.14} 
we can take a sufficiently small $\epsilon>0$ so that 
\begin{equation}
\label{eq:3.17}
\rho^{\frac{r}{p-1}-N}\Psi_{r,\rho}[v_\epsilon](T)\le 2\Lambda,
\end{equation}
where $v_\epsilon:=\max\{\pm v,0\}+\epsilon$. 
Similarly to \eqref{eq:3.8}, 
for any $0<\tau<t\le T$, 
multiplying \eqref{eq:1.1} by $v_\epsilon^{r-1}\zeta^k$ and integrating it in $\Omega\times(\tau,t)$, 
we obtain 
\begin{equation}
\label{eq:3.18}
\begin{split}
 & \int_{\Omega(x,2\rho)}v_\epsilon(y,s)^r\zeta^k\,dy\biggr|_{s=\tau}^{s=t}
+\int_\tau^t\int_{\Omega(x,\rho)}|\nabla v_\epsilon^{\frac{r}{2}}|^2\,dyds\\
 & \le C\rho^{-2}\int_\tau^t\int_{\Omega(x,2\rho)} v_\epsilon^r\,dyds
 +C\int_\tau^t\int_{\partial\Omega(x,2\rho)}v_\epsilon^{p+r-1}\zeta^k\,d\sigma ds. 
\end{split}
\end{equation}
This together with $v\in C(\overline{\Omega}\times[\tau,T])\cap L^\infty(\tau,T:L^\infty(\Omega))$ 
(see Definition~\ref{Definition:1.1}) implies that 
\begin{equation}
\label{eq:3.19}
\sup_{x\in\overline{\Omega}}\,
\int_\tau^t\int_{\Omega(x,\rho)}|\nabla v_\epsilon^{\frac{r}{2}}|^2\,dyds<\infty.
\end{equation}
Furthermore, by Lemma~\ref{Lemma:2.4}, \eqref{eq:3.17} and \eqref{eq:3.18}
we have 
\begin{equation}
\label{eq:3.20}
\begin{split}
 & \int_{\Omega(x,2\rho)}v_\epsilon(y,s)^r\zeta^k\,dy\biggr|_{s=\tau}^{s=t}
+\int_\tau^t\int_{\Omega(x,\rho)}|\nabla v_\epsilon^{\frac{r}{2}}|^2\,dyds
\le C\rho^{-2}\int_\tau^t\int_{\Omega(x,2\rho)} v_\epsilon^r\,dyds\\
 & \qquad\qquad\qquad
 +C(2\Lambda)^{\frac{p-1}{r}}
\left[\sup_{x\in\overline{\Omega}}\,\int_\tau^t\int_{\Omega(x,\rho)}|\nabla v_\epsilon^{\frac{r}{2}}|^2\,dyds
+\rho^{-2}(t-\tau)\Psi_{r,\rho}[v_\epsilon](t)\right]
\end{split}
\end{equation}
for $0<\tau<t\le T$. 
Therefore, by Lemma~\ref{Lemma:2.1}, \eqref{eq:1.18} and \eqref{eq:3.20}
we obtain
\begin{equation}
\begin{split}
 & \sup_{x\in\overline{\Omega}}\,\int_{\Omega(x,2\rho)}v_\epsilon(y,t)^r\,dy
+\sup_{x\in\overline{\Omega}}\,\int_\tau^t\int_{\Omega(x,\rho)}|\nabla v_\epsilon^{\frac{r}{2}}|^2\,dyds\\
 & \le M\sup_{x\in\overline{\Omega}}\,\int_{\Omega(x,\rho)}v_\epsilon(y,\tau)^r\,dy
 +C\rho^{-2}(t-\tau)\Psi_{r,\rho}[v_\epsilon](t)\\
 &\qquad
+C(2\Lambda)^{\frac{p-1}{r}}
\left[\sup_{x\in\overline{\Omega}}\,\int_\tau^t\int_{\Omega(x,\rho)}|\nabla v_\epsilon^{\frac{r}{2}}|^2\,dyds
+\rho^{-2}(t-\tau)\Psi_{r,\rho}[v_\epsilon](t)\right]
\end{split}
\label{eq:3.21}
\end{equation}
for $0<\tau<t\le T$. 
Taking a sufficiently small $\Lambda$ if necessary, 
we deduce from \eqref{eq:3.19} and \eqref{eq:3.21} that
\begin{equation*}
\begin{split}
 & \sup_{x\in\overline{\Omega}}\,\int_{\Omega(x,\rho)}v_\epsilon(y,t)^r\,dy
 +\frac{1}{2}\sup_{x\in\overline{\Omega}}\,\int_\tau^t\int_{\Omega(x,\rho)}|\nabla v_\epsilon^{\frac{r}{2}}|^2\,dyds\\
 & \qquad
 \le M\sup_{x\in\overline{\Omega}}\,\int_{\Omega(x,\rho)}v_\epsilon(y,\tau)^r\,dy
 +C\rho^{-2}(t-\tau)\Psi_{r,\rho}[v_\epsilon](t).
\end{split}
\end{equation*} 
Taking a sufficiently small $\mu\in(0,1]$, 
we obtain 
\begin{equation}
\label{eq:3.22}
\begin{split}
 & \Psi_{r,\rho}[v_\epsilon](t)+\frac{1}{2}\sup_{x\in\overline{\Omega}}\,\int_\tau^t\int_{\Omega(x,\rho)}|\nabla v_\epsilon^{\frac{r}{2}}|^2\,dyds\\
 & \le 2M\Psi_{r,\rho}[v_\epsilon](\tau)+C\rho^{-2}(t-\tau)\Psi_{r,\rho}[v_\epsilon](t)
\le 2M\Psi_{r,\rho}[v_\epsilon](\tau)+\frac{1}{2}\Psi_{r,\rho}[v_\epsilon](t)
\end{split}
\end{equation}
for $0<\tau<t\le T$ with $t-\tau\le\mu\rho^2$. 
This implies that 
\begin{equation}
\label{eq:3.23}
\Psi_{r,\rho}[\max\{\pm v,0\}](t)\le\Psi_{r,\rho}[v_\epsilon](t)\le 4M\Psi_{r,\rho}[v_\epsilon](\tau)
\le 5M\Psi_{r,\rho}[v](\tau)+C\epsilon^r\rho^N
\end{equation}
for $0<\tau<t\le T$ with $t-\tau\le\mu\rho^2$. 
Furthermore, by Lemma~\ref{Lemma:2.4}, \eqref{eq:3.22} and \eqref{eq:3.23}
we have 
\begin{equation}
\label{eq:3.24}
\begin{split}
 & \int_\tau^t\int_{\partial\Omega(x,\rho)}\max\{\pm v,0\}^{p+r-1}\,d\sigma ds
 \le \int_\tau^t\int_{\partial\Omega(x,\rho)}v_\epsilon^{p+r-1}\,d\sigma ds\\
 & \qquad\quad
\le C\Lambda^{\frac{p-1}{r}}\Psi_{r,\rho}[v_\epsilon](\tau)
\le C\Lambda^{\frac{p-1}{r}}\Psi_{r,\rho}[v](\tau)+C\epsilon^r\rho^N.
\end{split}
\end{equation} 
Since $\tau$ and $\epsilon$ is arbitrary, 
by \eqref{eq:3.23} and \eqref{eq:3.24} 
we obtain \eqref{eq:3.15} and \eqref{eq:3.16}. Thus Lemma~\ref{Lemma:3.2} follows. 
$\Box$
\begin{lemma}
\label{Lemma:3.3}
Assume the same conditions as in Lemma~{\rm\ref{Lemma:3.1}}. 
Let $r$ satisfy \eqref{eq:1.7} and $r>1$. 
Then there exists a positive constant $\Lambda$ such that, 
if 
\begin{equation}
\label{eq:3.25}
\rho^{\frac{r}{p-1}-N}\left(\Psi_{r,\rho}[v](T)+\Psi_{r,\rho}[w](T)\right)\le\Lambda
\end{equation}
for some $\rho\in(0,\rho_*/2)$, 
then 
\begin{equation}
\label{eq:3.26}
\Psi_{r,\rho}[z_0](t)\le C\Psi_{r,\rho}[z_0](\tau)
\end{equation}
for $0\le\tau<t\le T$ with $t-\tau\le\mu\rho^2$, 
where $C$ and $\mu$ are positive constants depending only on $N$, $\Omega$, $p$ and $r$.  
\end{lemma}
{\bf Proof.}
Let $x\in\overline{\Omega}$ and $\zeta$ be as in Lemma~\ref{Lemma:2.4}. 
Let $k$ be as in Lemma~\ref{Lemma:2.4} and $\epsilon>0$. 
Similarly to \eqref{eq:3.18}, 
we have
\begin{equation}
\label{eq:3.27}
\begin{split}
 & \int_{\Omega(x,2\rho)} z_\epsilon(y,s)^r\zeta^k\,dy\biggr|_{s=\tau}^{s=t}
+\int_\tau^t\int_{\Omega(x,2\rho)}|\nabla z_\epsilon^{\frac{r}{2}}|^2\zeta^k\,dyds\\
 & \le  C\rho^{-2}\int_\tau^t\int_{\Omega(x,2\rho)}z_\epsilon^r\,dyds
+C\int_\tau^t\int_{\partial\Omega(x,2\rho)}a(y,s)z_\epsilon^r\zeta^k\,d\sigma ds
\end{split}
\end{equation}
for all $0<\tau<t\le T$. 
This together with $z_\epsilon$, $a\in C(\overline{\Omega}\times[\tau,T])\cap L^\infty(\Omega\times(\tau,T))$ implies that 
\begin{equation}
\label{eq:3.28}
\sup_{x\in\overline{\Omega}}\,\int_\tau^t\int_{\Omega(x,2\rho)}|\nabla z_\epsilon^{\frac{r}{2}}|^2\,dyds<\infty
\end{equation}
for $0<\tau<t\le T$. 
On the other hand, 
by the H\"older inequality and \eqref{eq:3.3} 
we have 
\begin{equation}
\label{eq:3.29}
\begin{split}
\int_\tau^t\int_{\partial\Omega(x,2\rho)}a(y,\tau)z_\epsilon^r\zeta^k\,d\sigma ds
 & \le C\left(\int_\tau^t\int_{\partial\Omega(x,2\rho)}(|v|^{p+r-1}+|w|^{p+r-1})\,d\sigma ds\right)^{\frac{p-1}{p+r-1}}\\
 & \qquad\qquad
\times\left(\int_\tau^t\int_{\partial\Omega(x,2\rho)}z_\epsilon^{p+r-1}\zeta^k\,d\sigma ds\right)^{\frac{r}{p+r-1}}. 
\end{split}
\end{equation}
Let $\Lambda$ and $\mu$ be sufficiently small positive constants. 
Then, by Lemma~\ref{Lemma:2.1}, \eqref{eq:3.16} and \eqref{eq:3.25} 
we see that
\begin{equation}
\label{eq:3.30}
\begin{split}
 & \int_\tau^t\int_{\partial\Omega(x,2\rho)}(|v|^{p+r-1}+|w|^{p+r-1})\,d\sigma ds\\
 & \le M\sup_{x\in\overline{\Omega}}\, \int_\tau^t\int_{\partial\Omega(x,\rho)}(|v|^{p+r-1}+|w|^{p+r-1})\,d\sigma ds\\
 & \le C\Lambda^{\frac{p-1}{r}}\left\{\Psi_{r,\rho}[v](\tau)+\Psi_{r,\rho}[w](\tau)\right\}
 \le C\Lambda^{\frac{p+r-1}{r}}\rho^{-\frac{r}{p-1}+N}
\end{split}
\end{equation}
for all $0<\tau<t\le T$ with $t-\tau\le\mu\rho^2$. 
Similarly, by Lemma~\ref{Lemma:2.4} we obtain 
\begin{equation}
\label{eq:3.31}
\begin{split}
 & \int_\tau^t\int_{\partial\Omega(x,2\rho)}z_\epsilon^{p+r-1}\zeta^k\,d\sigma ds
 \le C\left(\rho^{\frac{r}{p-1}-N}\Psi_{r,\rho}[z_\epsilon](t)\right)^{\frac{p-1}{r}}\\
 & \qquad\qquad\times
 \biggr[\sup_{x\in\overline{\Omega}}\,\int_\tau^t\int_{\Omega(x,\rho)}|\nabla(z_\epsilon)^{\frac{r}{2}}|^2\,dyds
 +\rho^{-2}(t-\tau)\Psi_{r,\rho}[z_\epsilon](\tau)\biggr]
\end{split}
\end{equation}
for all $0<\tau<t\le T$ with $t-\tau\le\mu\rho^2$.  
Then we deduce from \eqref{eq:3.29}--\eqref{eq:3.31} that
\begin{equation}
\label{eq:3.32}
\begin{split}
 & \int_\tau^t\int_{\partial\Omega(x,2\rho)}a(y,t)z_\epsilon^r\zeta^k\,d\sigma ds\\
 & \le C\Lambda^{\frac{p-1}{r}}\left(\Psi_{r,\rho}[z_\epsilon](t)\right)^{\frac{p-1}{p+r-1}}\\
 & \qquad\times
 \biggr[\sup_{x\in\overline{\Omega}}\int_\tau^t\int_{\Omega(x,\rho)}|\nabla(z_\epsilon)^{\frac{r}{2}}|^2\,dyds
 +\rho^{-2}(t-\tau)\Psi_{r,\rho}[z_\epsilon](t)\biggr]^{\frac{r}{p+r-1}}\\
 & \le C\Lambda^{\frac{p-1}{r}}
  \biggr[\sup_{x\in\overline{\Omega}}\int_\tau^t\int_{\Omega(x,\rho)}|\nabla z_\epsilon^{\frac{r}{2}}|^2\,dyds
  +\Psi_{r,\rho}[z_\epsilon](t)+\rho^{-2}(t-\tau)\Psi_{r,\rho}[z_\epsilon](\tau)\biggr]
\end{split}
\end{equation}
for all $0<\tau<t\le T$ with $t-\tau\le\mu\rho^2$. 
Therefore, by Lemma~\ref{Lemma:2.1}, \eqref{eq:3.27} and \eqref{eq:3.32} 
we have 
\begin{equation*}
\begin{split}
  & \sup_{x\in\overline{\Omega}}\,\int_{\Omega(x,\rho)}z_\epsilon^r\,dy
 +\sup_{x\in\overline{\Omega}}\,\int_\tau^t\int_{\Omega(x,\rho)}|\nabla z_\epsilon^{\frac{r}{2}}|^2\,dyds\\
 & \quad
 \le M\Psi_{r,\rho}[z_\epsilon](\tau)+C\rho^{-2}(t-\tau)\Psi_{r,\rho}[z_\epsilon](t)\\
 & \qquad
 +C\Lambda^{\frac{p-1}{r}}
 \biggr[\sup_{x\in\overline{\Omega}}\int_\tau^t\int_{\Omega(x,\rho)}|\nabla z_\epsilon^{\frac{r}{2}}|^2\,dyds
 +\Psi_{r,\rho}[z_\epsilon](t)+\rho^{-2}(t-\tau)\Psi_{r,\rho}[z_\epsilon](\tau)\biggr]
\end{split}
\end{equation*}
for all $0<\tau<t\le T$ with $t-\tau\le\mu\rho^2$. 
Then, taking sufficiently small constants $\Lambda$ and $\mu$ if necessary, 
we obtain 
$$
\Psi_{r,\rho}[z_\epsilon](t)\le 4M\Psi_{r,\rho}[z_\epsilon](\tau)
$$
for all $0<\tau<t\le T$ with $t-\tau\le\mu\rho^2$. 
This implies \eqref{eq:3.26}, and the proof is complete. 
$\Box$\vspace{5pt}

Now we are ready to complete the proof of Theorems~\ref{Theorem:1.1} and \ref{Theorem:1.2} in the case $r>1$. 
\vspace{5pt}
\newline
{\bf Proof of Theorem~\ref{Theorem:1.1} in the case $r>1$.} 
Let $\gamma_1$ be a sufficiently small positive constant and assume \eqref{eq:1.8}. 
Let $\{\varphi_n\}$ satisfy \eqref{eq:1.14} and define $T_n^*$ and $T_n^{**}$ as in \eqref{eq:1.18}. 
Then it follows from \eqref{eq:1.17} that 
\begin{equation}
\label{eq:3.33}
\rho^{\frac{r}{p-1}-N}\Psi_{r,\rho}[u_n](t)\le 6M\rho^{\frac{r}{p-1}-N}\Psi_{r,\rho}[u_n](0)\le 6M(2\gamma_1)^r
\end{equation}
for all $0\le t\le T_n^*$. 
Taking a sufficiently small $\gamma_1$ if necessary, 
by Lemma~\ref{Lemma:3.2}, \eqref{eq:1.17}  and \eqref{eq:3.33}, 
we can find a constant $\mu>0$ such that 
\begin{equation}
\label{eq:3.34}
\Psi_{r,\rho}[u_n](t)\le 5M\Psi_{r,\rho}[u_n](0)
<6M\Psi_{r,\rho}[u_n](0)\le C\|\varphi\|_{r,\rho}^r
\end{equation}
for $0\le t\le\min\{T_n^*,\mu\rho^2\}$. 
On the other hand, 
we apply Lemma~\ref{Lemma:3.1} with $R_1=\rho/2$, $R_2=\rho$, $t_1=t/2$ and $t_2=t/4$ 
to obtain 
\begin{eqnarray}
\label{eq:3.35}
 & & \|u_n(t)\|_{L^\infty(\Omega(x,\rho/2))}\le CD^{\frac{N+2}{2r}}
\left(\int_{t/4}^t\int_{\Omega(x,\rho)}|u_n|^r\,dyds\right)^{1/r},\\
\label{eq:3.36}
 & & \int_{t/2}^t\int_{\Omega(x,\rho/2)}|\nabla u_n|^2\,dyds
 \le CD\int_{t/4}^t\int_{\Omega(x,\rho)}|u_n|^2\,dyds,
\end{eqnarray}
for all $x\in\overline{\Omega}$ and $t\in(0,T_n)$. 
where 
$D=\||u_n|^{p-1}\|_{L^\infty(\Omega(x,\rho)\times(t/4,t))}^2+\rho^{-2}+t^{-1}$. 
By \eqref{eq:1.18}, \eqref{eq:3.34} and \eqref{eq:3.35} 
we have 
\begin{eqnarray}
\label{eq:3.37}
 & & \|u_n(t)\|_{L^\infty(\Omega)}
\le Ct^{-\frac{N}{2r}}\|\varphi\|_{r,\rho}
\le C\gamma_1t^{-\frac{1}{2(p-1)}}(\rho^{-2}t)^{-\frac{N}{2r}+\frac{1}{2(p-1)}},\\
\label{eq:3.38}
 & & \sup_{x\in\overline{\Omega}}\int_{t/2}^t\int_{\Omega(x,\rho)}
|\nabla u_n|^2\,dyds
\le C\rho^N\|u_n\|_{L^\infty(\Omega\times(t/4,t))}^2
\le C\rho^Nt^{-\frac{N}{r}}\|\varphi\|_{r,\rho}^2,\qquad
\end{eqnarray}
for all $0<t\le\min\{\mu\rho^2,T_n^*,T_n^{**}\}$.  
Since $r\ge N(p-1)$, 
taking sufficiently small $\gamma_1>0$ and $\mu>0$ if necessary, 
by \eqref{eq:3.37} we have 
$$
(\rho^{-2}t)^{\frac{1}{2}}+t^{\frac{1}{2}}\|u_n(t)\|_{L^\infty(\Omega)}^{p-1}
\le\mu^{\frac{1}{2}}+(C\gamma_1)^{p-1}\mu^{-\frac{N(p-1)}{2r}+\frac{1}{2}}\le 1
$$
for $0<t\le\min\{\mu\rho^2,T_n^*,T_n^{**}\}$. 
This implies that 
$T_n>T_n^{**}>\min\{T_n^*,\mu\rho^2\}$ for $n=1,2,\dots$.
Then, by \eqref{eq:3.34} we see that 
$T_n^*>\mu\rho^2$ for $n=1,2,\dots$. 
Therefore, by \eqref{eq:3.34}, \eqref{eq:3.37} and \eqref{eq:3.38}
we obtain 
\begin{eqnarray}
\label{eq:3.39}
 & & \|u_n(t)\|_{L^\infty(\Omega)}\le Ct^{-\frac{N}{2r}}\|\varphi\|_{r,\rho},\\
\label{eq:3.40}
 &  & \sup_{x\in\overline{\Omega}}\int_{t/2}^t\int_{\Omega(x,\rho)}
|\nabla u_n|^2\,dyds\le C\rho^Nt^{-\frac{N}{r}}\|\varphi\|_{r,\rho}^2,\\
\label{eq:3.41}
 & & \sup_{0<t<\mu\rho^2}\|u_n(t)\|_{r,\rho}\le C\|\varphi\|_{r,\rho},
\end{eqnarray}
for $0<t\le\mu\rho^2$ and $n=1,2,\dots$. 

Applying \cite[Theorem~6.2]{DB01} with the aid of \eqref{eq:3.39}, 
we see that 
$u_n$ $(n=1,2,\dots)$ are uniformly bounded and equicontinuous on $K\times[\tau,\mu\rho^2]$ 
for any compact set $K\subset\overline{\Omega}$ and $\tau\in(0,\mu\rho^2]$. 
Then, by the Ascoli-Arzel\`a theorem and the diagonal argument
we can find a subsequence $\{u_{n'}\}$ and a continuous function $u$ in $\Omega\times(0,\mu\rho^2]$ 
such that 
$$
\lim_{n'\to\infty}\,\|u_{n'}-u\|_{L^\infty(K\times[\tau,\mu\rho^2])}=0
$$
for any compact set $K\subset\overline{\Omega}$ and $\tau\in(0,\mu\rho^2]$. 
This together with \eqref{eq:3.39} and \eqref{eq:3.41} implies 
\eqref{eq:1.9} and \eqref{eq:1.10}. 
Furthermore, by  \eqref{eq:3.40},
taking a subsequence if necessary, 
we see that 
$$
\lim_{n'\to\infty}u_{n'}=u\quad\mbox{weakly in}\,\,\,L^2([\tau,\mu\rho^2]:W^{1,2}(\Omega\cap B(0,R)))
$$
for any $R>0$ and $0<\tau<\mu\rho^2$. This implies that $u$ satisfies \eqref{eq:1.6}. 

On the other hand, since $u_n$ is a $L^r_{uloc}(\Omega)$-solution of \eqref{eq:1.1} (see \eqref{eq:1.16}), 
we see that
$$
u_n\in C([0,\mu\rho^2]:L^r_{uloc,\rho}(\Omega)). 
$$
Furthermore, 
by Lemma~\ref{Lemma:3.3} and \eqref{eq:3.33}, 
taking a sufficiently small $\gamma_1$ if necessary, 
we have  
$$
\sup_{0<\tau<\mu\rho^2}\|u_m(\tau)-u_n(\tau)\|_{r,\rho}\le C\|u_m(0)-u_n(0)\|_{r,\rho},
\quad m,n=1,2,\dots. 
$$
This means that 
$\{u_n\}$ is a Cauchy sequence in $C([0,\mu\rho^2]:L^r_{uloc,\rho}(\Omega))$, 
which implies 
\begin{equation}
\label{eq:3.42}
u\in C([0,\mu\rho^2]:L^r_{uloc,\rho}(\Omega)). 
\end{equation}
Therefore we see that 
$u$ is a $L^r_{uloc}(\Omega)$-solution of \eqref{eq:1.1} in $\Omega\times[0,\mu\rho^2]$ 
satisfying \eqref{eq:1.9} and \eqref{eq:1.10}, 
and the proof of Theorem~\ref{Theorem:1.1} for the case $r>1$ is complete. 
$\Box$
\vspace{5pt}

\noindent
{\bf Proof of Theorem~\ref{Theorem:1.2} in the case $r>1$.} 
Let $v$ and $w$ be $L^r_{uloc}(\Omega)$-solutions of \eqref{eq:1.1} 
in $\Omega\times[0,T)$, where $T>0$. 
Let $\gamma_2$ be a sufficiently small constant and assume \eqref{eq:1.12}. 
We can assume, without loss of generality, that $\rho\in(0,\rho_*/2)$. 
Since $v$, $w\in C([0,T]:L^r_{uloc,\rho}(\Omega))$, 
we can find a constant $T'\in(0,T)$ such that 
\begin{equation}
\label{eq:3.43}
\rho^{\frac{1}{p-1}-\frac{N}{r}}\left[\sup_{0<\tau\le T'}\|v(\tau)\|_{r,\rho}+\sup_{0<\tau\le T'}\|w(\tau)\|_{r,\rho}\right]\le 2\gamma_2.
\end{equation}
Furthermore, for any $T''\in(T',T)$, since $v$, $w\in L^\infty(\Omega\times(T',T''))$, 
we see that  
\begin{equation}
\label{eq:3.44}
\tilde{\rho}^{\frac{1}{p-1}-\frac{N}{r}}
\left[\sup_{T'<\tau\le T''}\|v(\tau)\|_{r,\tilde{\rho}}+\sup_{T_1<\tau\le T_2}\|w(\tau)\|_{r,\tilde{\rho}}\right]\le\gamma_2
\end{equation}
for some $\tilde{\rho}\in(0,\rho)$. 
Since $v(x,0)\le w(x,0)$ for almost all $x\in\Omega$, 
by \eqref{eq:3.43} and \eqref{eq:3.44} 
we apply Lemma~\ref{Lemma:3.3} to obtain 
$$
\sup_{0<\tau<\min\{\mu\tilde{\rho}^2,T''\}}\|(v(\tau)-w(\tau))_+\|_{r,\tilde{\rho}}\le C\|(v(0)-w(0))_+\|_{r,\tilde{\rho}}=0
$$
for some constant $\mu>0$. This implies that $v(x,t)\le w(x,t)$ in $\Omega\times(0,\min\{\mu\tilde{\rho}^2,T''\}]$. 
Repeating this argument, we see that $v(x,t)\le w(x,t)$ in $\Omega\times(0,T'']$. 
Finally, since $T''$ is arbitrary, we see that $v(x,t)\le w(x,t)$ in $\Omega\times(0,T)$, 
and the proof is complete. 
$\Box$
\section{Proof of Theorems~\ref{Theorem:1.1} and \ref{Theorem:1.2} in the case $r=1$}
In this section we consider the case $1<p<1+1/N$ and $r=1$, 
and complete the proof of Theorems~\ref{Theorem:1.1} and \ref{Theorem:1.2}. 
Furthermore, we prove Corollary~\ref{Corollary:1.1}. 
We use the same notation as in Section~3. 
\begin{lemma}
\label{Lemma:4.1}
Assume the same conditions as in Theorem~{\rm\ref{Theorem:1.1}}. 
Let $v$ and $w$ be $L^1_{uloc}(\Omega)$-solutions 
of \eqref{eq:1.1} in $\Omega\times[0,T]$, where $0<T<\infty$, such that 
\begin{equation}
\label{eq:4.1}
\|v(t)\|_{L^\infty(\Omega)}+\|w(t)\|_{L^\infty(\Omega)}
\le C_1t^{-\frac{1}{2(p-1)}},\qquad 0<t\le T,
\end{equation}
for some $C_1>0$. 
Then there exists a constant $C_2$ such that  
\begin{eqnarray}
\label{eq:4.2}
 & & \|v(t)\|_{L^\infty(\Omega)}\le C_2t^{-\frac{N}{2}}\Psi_{1,\rho}[v](t),\\
\label{eq:4.3}
 & & \|z_0(t)\|_{L^\infty(\Omega)}\le C_2t^{-\frac{N}{2}}\Psi_{1,\rho}[z_0](t),
\end{eqnarray}
for all $0<t\le\min\{T,\rho^2\}$ and $0<\rho<\rho_*$. 
\end{lemma}
{\bf Proof.}
Similarly to \eqref{eq:3.35}, 
by Lemma~\ref{Lemma:3.1} and \eqref{eq:4.1} 
we have 
\begin{equation*}
\begin{split}
\|z_0(t)\|_{L^\infty(\Omega(x,\rho/2))}
 & \le C\left[\|v(t)\|_{L^\infty(\Omega\times(t/4,t))}^{2(p-1)}+\rho^{-2}+t^{-1}\right]^{\frac{N+2}{2}}
\int_{t/4}^t\int_{\Omega(x,\rho)} |z_0(y,s)|\,dyds\\
 & \le C(1+C_1^{2(p-1)})^{\frac{N+2}{2}}t^{-\frac{N}{2}}\Psi_{1,\rho}[z_0](t)
\end{split}
\end{equation*}
for all $x\in\overline{\Omega}$ and $0<t\le\min\{T,\rho^2\}$. 
This implies \eqref{eq:4.3}. 
Furthermore, \eqref{eq:4.2} follows from \eqref{eq:4.3}, 
and the proof is complete.
$\Box$
\begin{lemma}
\label{Lemma:4.2}
Assume the same conditions as in Theorem~{\rm\ref{Theorem:1.1}} and $1<p<1+1/N$. 
Let $v$ and $w$ be $L^1_{uloc}(\Omega)$-solutions 
of \eqref{eq:1.1} in $\Omega\times[0,T]$, where $0<T<\infty$, 
and assume \eqref{eq:4.1} for some constant $C_1>0$. 
Let $0<\rho<\rho_*$ and $\Lambda$ be such that 
\begin{equation}
\label{eq:4.4}
\rho^{\frac{1}{p-1}-N}\left[\Psi_{1,\rho}[v](T)+\Psi_{1,\rho}[w](T)\right]\le\Lambda. 
\end{equation}
Then, for any $\sigma\in(0,1)$ and $\delta\in(0,1)$ with $\sigma>\delta N/2$, 
there exists a positive constant $C_2$ such that 
\begin{equation}
\label{eq:4.5}
\limsup_{\epsilon\to 0}\,\sup_{x\in\overline{\Omega}}\,\int_0^t \int_{\Omega(x,\rho)}  
(\rho^{-2}s)^\sigma\frac{|\nabla z_\epsilon|^2}{z_\epsilon^{1-\delta}}\,dyds
\le C_2\mu^{\sigma-\frac{\delta N}{2}}\rho^{-\delta N}\Psi_{1,\rho}[z_0](t)^{1+\delta}
\end{equation}
for $0<t\le\min\{T,\mu\rho^2\}$ and $0<\mu\le 1$. 
\end{lemma}
{\bf Proof.}
Let $\sigma\in(0,1)$ and $\delta\in(0,1)$ be such that $\sigma>\delta N/2$. 
Let $x\in\overline{\Omega}$ and $\zeta$ be as in Lemma~{\rm\ref{Lemma:2.4}}. 
Similarly to \eqref{eq:3.8}, 
multiplying \eqref{eq:3.1} by $(\rho^{-2}t)^\sigma z_\epsilon(x,t)^\delta\zeta(x)^2$ 
and integrating it on $\Omega(x,2\rho)\times(\tau,t)$, 
we obtain 
\begin{equation}
\begin{split}
 & \frac{\delta}{2} \int_\tau^t \int_{\Omega(x,2\rho)}  
(\rho^{-2}s)^\sigma\frac{|\nabla z_\epsilon|^2}{z_\epsilon^{1-\delta}}\zeta^2\,dyds
\le \frac{(\rho^{-2}\tau)^\sigma}{1+\delta}\int_{\Omega(x,2\rho)}z_\epsilon(y,\tau)^{1+\delta}\,dy\\
 & \qquad
 +\frac{\sigma}{1+\delta}\rho^{-2}\int_\tau^t \int_{\Omega(x,2\rho)}(\rho^{-2}s)^{\sigma-1} z_\epsilon^{1+\delta} \zeta^2 \,dyds\\
 & \qquad
 +C\rho^{-2}\int_\tau^t \int_{\Omega(x,2\rho)} 
 (\rho^{-2}s)^\sigma z_\epsilon^{1+\delta} \,dyds
 +\int_\tau^t\int_{\partial\Omega(x,2\rho)}(\rho^{-2}s)^\sigma a(y,s)z_\epsilon^{1+\delta}\zeta^2 \,d\sigma
\end{split}
\label{eq:4.6}
\end{equation}
for $0<\tau<t\le T$. 
On the other hand, 
it follows from Lemma~\ref{Lemma:4.1}, \eqref{eq:3.3}, \eqref{eq:4.1} and  \eqref{eq:4.4} that 
\begin{equation}
\label{eq:4.7}
\|a(t)\|_{L^\infty(\Omega)}\le Ct^{-\frac{N(p-1)}{2}}\left[\Psi_{1,\rho}[v](t)^{p-1}+\Psi_{1,\rho}[w](t)^{p-1}\right]
\le C\Lambda^{p-1}\rho^{-1}(\rho^{-2}t)^{-\frac{N(p-1)}{2}}
\end{equation}
for all $0<t\le\min\{T,\rho^2\}$. 
Furthermore, 
by Lemma~\ref{Lemma:2.3} we have 
\begin{equation}
\label{eq:4.8}
\begin{split}
 & \int_{\partial\Omega(x,2\rho)}z_\epsilon^{1+\delta}\zeta^2 \,d\sigma
 \le\nu\int_{\Omega(x,2\rho)}|\nabla(z_\epsilon^{\frac{1+\delta}{2}}\zeta)|^2 \,dy
+\frac{C}{\nu}\int_{\Omega(x,2\rho)}(z_\epsilon^{\frac{1+\delta}{2}}\zeta)^2 \,dy\\
 & \qquad\qquad
 \le 2\nu\left(\frac{1+\delta}{2}\right)^2\int_{\Omega(x,2\rho)}\frac{|\nabla z_\epsilon|^2}{z_\epsilon^{1-\delta}}\zeta^2\,dy
 +2\nu\int_{\Omega(x,2\rho)}z_\epsilon^{1+\delta}|\nabla\zeta|^2\,dy\\
 & \qquad\qquad\qquad\qquad\qquad\qquad\qquad\qquad\qquad\quad\,\,\,
+\frac{C}{\nu}\int_{\Omega(x,2\rho)}z_\epsilon^{1+\delta}\zeta^2 \,dy
\end{split}
\end{equation}
for all $0<t\le T$ and $\nu>0$. 
By \eqref{eq:4.7} and \eqref{eq:4.8} we obtain
\begin{equation}
\label{eq:4.9}
\begin{split}
 & \int_\tau^t\int_{\partial\Omega(x,2\rho)}(\rho^{-2}s)^\sigma a(y,s)z_\epsilon^{1+\delta}\zeta^2 \,d\sigma ds\\
 & \le\frac{\delta}{4}\int_\tau^t \int_{\Omega(x,2\rho)}  
 (\rho^{-2}s)^\sigma\frac{|\nabla z_\epsilon|^2}{z_\epsilon^{1-\delta}}\zeta^2\,dyds
 +C\rho^{-2}\int_\tau^t \int_{\Omega(x,2\rho)}(\rho^{-2}s)^\sigma z_\epsilon^{1+\delta}\,dyds\\
  & \qquad\qquad\qquad\qquad\qquad
 +C\rho^{-2}\int_\tau^t \int_{\Omega(x,2\rho)} 
 (\rho^{-2}s)^{\sigma-N(p-1)}z_\epsilon^{1+\delta}\,dyds
\end{split}
\end{equation}
for all $0<\tau<t\le\min\{T,\rho^2\}$. 
We deduce from \eqref{eq:4.6}--\eqref{eq:4.9} that 
\begin{equation}
\label{eq:4.10}
\begin{split}
 & \frac{\delta}{4}\int_\tau^t \int_{\Omega(x,2\rho)}  
(\rho^{-2}s)^\sigma\frac{|\nabla z_\epsilon|^2}{z_\epsilon^{1-\delta}}\zeta^2\,dyds
\le \frac{(\rho^{-2}\tau)^\sigma}{1+\delta}\int_{\Omega(x,2\rho)}z_\epsilon(y,\tau)^{1+\delta}\,dy\\
 & \qquad
+C\rho^{-2}\int_\tau^t \int_{\Omega(x,2\rho)} 
[(\rho^{-2}s)^{\sigma-1} +(\rho^{-2}s)^\sigma+ (\rho^{-2}s)^{\sigma-N(p-1)}]
z_\epsilon^{1+\delta}\,dyds
\end{split}
\end{equation}
for all $0<\tau<t\le\min\{T,\rho^2\}$. 
Furthermore, 
by Lemmas~\ref{Lemma:2.1} and \ref{Lemma:4.1} 
we have  
\begin{equation}
\label{eq:4.11}
\begin{split}
 & \sup_{x\in\overline{\Omega}}\,\int_{\Omega(x,2\rho)}z_\epsilon(y,s)^{1+\delta}\,dy
 \le 2M\sup_{x\in\overline{\Omega}}\,\int_{\Omega(x,\rho)}z_0(y,s)^{1+\delta}\,dy+C\epsilon^{1+\delta}\rho^N\\
 & \qquad
 \le 2M\|z_0(s)\|_{L^\infty(\Omega)}^\delta\Psi_{1,\rho}[z_0](t)+C\epsilon^{1+\delta}\rho^N\\
 & \qquad
 \le C(\rho^{-2}s)^{-\frac{\delta N}{2}}\rho^{-\delta N}\Psi_{1,\rho}[z_0](t)^{1+\delta}+C\epsilon^{1+\delta}\rho^N
\end{split}
\end{equation}
for all $0<s<t\le\min\{T,\rho^2\}$. 
It follows from $N(p-1)<1$ and $\sigma>\delta N/2$ that 
$$
\sigma-N(p-1)-\frac{\delta N}{2}>\sigma-1-\frac{\delta N}{2}>-1.
$$
Then, by \eqref{eq:4.10} and \eqref{eq:4.11}, 
passing to the limit as $\tau\to 0$ and $\epsilon\to 0$, 
we have 
\begin{equation*}
\begin{split}
 & \limsup_{\epsilon\to 0}\,\sup_{x\in\overline{\Omega}}\,
\int_0^t \int_{\Omega(x,\rho)}  
(\rho^{-2}s)^\sigma\frac{|\nabla z_\epsilon|^2}{z_\epsilon^{1-\delta}}\,dyds\\
 & \le C\rho^{-2-\delta N}\Psi_{1,\rho}[z_0](t)^{1+\delta}\int_0^t (\rho^{-2}s)^{-\frac{\delta N}{2}}
[(\rho^{-2}s)^{\sigma-1} +(\rho^{-2}s)^\sigma+ (\rho^{-2}s)^{\sigma-N(p-1)}]\,ds\\
& \le C\rho^{-\delta N}\mu^{\sigma-\frac{\delta N}{2}}\Psi_{1,\rho}[z_0](t)^{1+\delta}
\end{split}
\end{equation*}
for all $0<t\le\min\{T,\mu\rho^2\}$ and $0<\mu\le 1$.  
This implies \eqref{eq:4.5}, and Lemma~\ref{Lemma:4.2} follows. 
$\Box$\vspace{5pt}
\begin{lemma}
\label{Lemma:4.3}
Assume the same conditions as in Lemma~{\rm\ref{Lemma:4.2}} with $\rho\in(0,\rho_*/2)$. 
Then there exists a constant $\mu\in(0,1)$ such that
\begin{equation}
\label{eq:4.12}
\Psi_{1,\rho}[z_0](t)\le 2M\Psi_{1,\rho}[z_0](0),
\quad 0<t\le\min\{T,\mu\rho^2\}.
\end{equation}
\end{lemma}
{\bf Proof.}
Let $x\in\overline{\Omega}$ and $\zeta$ be as in Lemma~\ref{Lemma:2.4}. 
Let $\sigma\in(0,1)$ and $\delta\in(0,1)$ be such that 
\begin{equation}
\label{eq:4.13}
\frac{\delta N}{2}<\sigma<1-N(p-1)
\qquad\mbox{and}\qquad
p-1>\delta.
\end{equation}
By \eqref{eq:3.1} we have  
\begin{equation}
\label{eq:4.14}
\int_{\Omega(x,2\rho)}z_0\zeta^2\,dy\biggr|_{s=\tau}^{s=t}
\le 2\int_\tau^t\int_{\Omega(x,2\rho)}|\nabla z_0||\nabla\zeta|\zeta\,dyds
+\int_\tau^t\int_{\partial\Omega(x,2\rho)}a(y,s)z_0\zeta^2\,d\sigma ds
\end{equation}
for $0<\tau<t\le T$. 
Furthermore, we have 
\begin{equation}
\label{eq:4.15}
\begin{split}
2\int_\tau^t\int_{\Omega(x,2\rho)}|\nabla z_0||\nabla\zeta|\zeta\,dyds
 & \le\nu\limsup_{\epsilon\to 0}\int_\tau^t\int_{\Omega(x,2\rho)}(\rho^{-2}s)^\sigma\frac{|\nabla z_\epsilon|^2}{z_\epsilon^{1-\delta}}\zeta^2\,dyds\\
 & \qquad\quad
+C\nu^{-1}\rho^{-2}\int_\tau^t\int_{\Omega(x,2\rho)}(\rho^{-2}s)^{-\sigma}z_0^{1-\delta}\,dyds
\end{split}
\end{equation}
for $\nu>0$. 
On the other hand, by \eqref{eq:1.3} and \eqref{eq:4.7} we obtain 
\begin{equation}
\begin{split}
 & \int_\tau^t\int_{\partial\Omega(x,2\rho)}a(y,s)z_0\zeta^2\,d\sigma ds
 \le\int_\tau^t\int_{\partial\Omega(x,2\rho)}a(y,s)z_\epsilon\zeta^2\,d\sigma ds\\
 & \le C\Lambda^{p-1}\rho^{-1}\int_\tau^t\int_{\partial\Omega(x,2\rho)}
 (\rho^{-2}s)^{-\frac{N(p-1)}{2}}z_\epsilon\zeta^2\,d\sigma ds\\
 & \le C\rho^{-1}\int_\tau^t(\rho^{-2}s)^{-\frac{N(p-1)}{2}}
 \int_{\Omega(x,2\rho)}[|\nabla z_\epsilon|\zeta^2+2z_\epsilon\zeta|\nabla\zeta|]\,dy ds\\
 & \le C\nu\int_\tau^t\int_{\Omega(x,2\rho)}
 (\rho^{-2}s)^\sigma\frac{|\nabla z_\epsilon|^2}{z_\epsilon^{1-\delta}}\zeta^2\,dyds\\
 & \qquad\qquad
 +C\rho^{-2}\nu^{-1}\int_\tau^t(\rho^{-2}s)^{-\sigma-N(p-1)}
 \int_{\Omega(x,2\rho)} z_\epsilon^{1-\delta}\,dyds\\
 & \qquad\qquad\qquad
 +C\rho^{-2}\int_\tau^t(\rho^{-2}s)^{-\frac{N(p-1)}{2}}
 \int_{\Omega(x,2\rho)} z_\epsilon\,dyds
\end{split}
\label{eq:4.16}
\end{equation}
for $0<t\le\min\{T,\rho^2\}$, $\epsilon>0$ and $\nu>0$. 
Then it follows from Lemma~\ref{Lemma:2.1} and \eqref{eq:4.14}--\eqref{eq:4.16} that  
\begin{equation}
\label{eq:4.17}
\begin{split}
 & \sup_{x\in\overline{\Omega}}\,\int_{\Omega(x,\rho)}z_0(y,t)\,dy\le M \sup_{x\in\overline{\Omega}}\,\int_{\Omega(x,\rho)}z_0(y,0)\,dy\\
 & \qquad
 +C\nu\limsup_{\epsilon\to 0}\,\sup_{x\in\overline{\Omega}}\,
 \int_0^t\int_{\Omega(x,2\rho)}(\rho^{-2}s)^\sigma
 \frac{|\nabla z_\epsilon|^2}{z_\epsilon^{1-\delta}}\zeta^2\,dyds\\
 & \qquad\qquad
 +C\nu^{-1}\rho^{-2}\sup_{x\in\overline{\Omega}}\,\int_0^t\int_{\Omega(x,\rho)}(\rho^{-2}s)^{-\sigma}z_0^{1-\delta}\,dyds\\
 & \qquad\qquad\qquad
 +C\rho^{-2}\nu^{-1}\sup_{x\in\overline{\Omega}}\,\int_0^t(\rho^{-2}s)^{-\sigma-N(p-1)}
 \int_{\Omega(x,\rho)} z_0^{1-\delta}\,dyds\\
 & \qquad\qquad\qquad\qquad
 +C\rho^{-2}\sup_{x\in\overline{\Omega}}\,\int_0^t(\rho^{-2}s)^{-\frac{N(p-1)}{2}}
 \int_{\Omega(x,\rho)} z_0\,dyds
\end{split}
\end{equation}
for $0<t\le\min\{T,\rho^2\}$ and $\nu>0$. 
Furthermore, by the H\"older inequality 
we have 
\begin{equation}
\label{eq:4.18}
\sup_{x\in\overline{\Omega}}\,
\int_{\Omega(x,\rho)} z_0(y,t)^{1-\delta}\,dy\le C\rho^{\delta N}\Psi_{1,\rho}[z_0](t)^{1-\delta},
\quad t>0.
\end{equation}
Then we deduce from \eqref{eq:4.5}, \eqref{eq:4.17} and \eqref{eq:4.18} that 
\begin{equation*}
\begin{split}
\Psi_{1,\rho}[z_0](t) & \le M\Psi_{1,\rho}[z_0](0)
+C\nu\mu^{\sigma-\frac{\delta N}{2}}\rho^{-\delta N}\Psi_{1,\rho}[z_0](t)^{1+\delta}\\
 & \quad
 +C\nu^{-1}\rho^{\delta N}\Psi_{1,\rho}[z_0](t)^{1-\delta}\rho^{-2}\int_0^t[(\rho^{-2}s)^{-\sigma}+(\rho^{-2}s)^{-\sigma-N(p-1)}]\,ds\\
 & \qquad
 +C\rho^{-2}\Psi_{1,\rho}[z_0](t)\int_0^t(\rho^{-2}s)^{-\frac{N(p-1)}{2}}\,ds
\end{split}
\end{equation*}
for $0<t\le\min\{T,\mu\rho^2\}$, $0<\mu\le 1$ and $\nu>0$. 
Then, taking $\nu=\rho^{\delta N}\Psi_{1,\rho}[z_0](t)^{-\delta}$ if $\Psi_{1,\rho}[z_0](t)\not=0$, 
we can find a positive constant $\mu\in(0,1)$ such that 
\begin{equation*}
\begin{split}
\Psi_{1,\rho}[z_0](t) & \le M\Psi_{1,\rho}[z_0](0)
+C(\mu^{\sigma-\frac{\delta N}{2}}+\mu^{1-\sigma-N(p-1)}+\mu^{1-\frac{N(p-1)}{2}})\Psi_{1,\rho}[z_0](t)\\
 & \le M\Psi_{1,\rho}[z_0](0)+\frac{1}{2}\Psi_{1,\rho}[z_0](t)
\end{split}
\end{equation*}
for $0<t\le\min\{T,\mu\rho^2\}$. 
This implies \eqref{eq:4.12}, and Lemma~\ref{Lemma:4.3} follows.
$\Box$\vspace{5pt}

Now we are ready to prove Theorem~\ref{Theorem:1.1} in the case $r=1$. 
\vspace{5pt}
\newline
{\bf Proof of Theorem~\ref{Theorem:1.1} in the case $r=1$.}
It suffices to consider the case $1<p<1+1/N$. 
Let $\gamma_1$ be a sufficiently small positive constant and assume \eqref{eq:1.8}. 
Let $\{\varphi_n\}$ satisfy \eqref{eq:1.14} and define $T_n^*$ and $T_n^{**}$ as in \eqref{eq:1.18}. 
Then it follows from \eqref{eq:1.17} that 
\begin{equation}
\label{eq:4.19}
\rho^{\frac{1}{p-1}-N}\Psi_{1,\rho}[u_n](t)\le 6M\rho^{\frac{1}{p-1}-N}\Psi_{1,\rho}[u_n](0)\le 12M\gamma_1
\end{equation}
for all $0\le t\le T_n^*$. 
By Lemma~\ref{Lemma:4.1} we have 
\begin{equation}
\label{eq:4.20}
\|u_n(t)\|_{L^\infty(\Omega)}\le Ct^{-\frac{N}{2}}\Psi_{1,\rho}[u_n](t)
\end{equation}
for $0<t\le\min\{T_n^{**},\rho^2\}<T_n$ and $n=1,2,\dots$.  
Then, taking a sufficiently small $\gamma_1$ and applying Lemma~\ref{Lemma:4.3} with $v=u_n$ and $w=0$, 
we can find a constant $\mu\in(0,1)$ such that  
\begin{equation}
\label{eq:4.21}
\Psi_{1,\rho}[u_n](t)\le 2M\Psi_{1,\rho}[u_n](0)
\end{equation}
for $0<t\le\min\{T_n^*,T_n^{**},\mu\rho^2\}$ and $n=1,2,\dots$. 
This implies that $\min\{T_n^{**},\mu\rho^2\}<T_n^*$ for $n=1,2,\dots$. 
Furthermore, by \eqref{eq:4.19}--\eqref{eq:4.21}, 
taking a sufficiently small $\mu$ if necessary, 
we obtain 
\begin{equation*}
\begin{split}
 (\rho^{-2}t)^{\frac{1}{2}}+t^{\frac{1}{2}}\|u_n(t)\|_{L^\infty(\Omega)}^{p-1}
 & \le \mu^{\frac{1}{2}}+C(\rho^{-2}t)^{-\frac{N(p-1)}{2}+\frac{1}{2}}\gamma_1^{p-1}\\
 & \le \mu^{\frac{1}{2}}+C\mu^{-\frac{N(p-1)}{2}+\frac{1}{2}}\gamma_1^{p-1}\le 1
\end{split}
\end{equation*}
for $0<t\le\min\{\mu\rho^2,T_n^{**}\}$. 
This yields $T_n^{**}>\mu\rho^2$ for $n=1,2,\dots$. 
Therefore, by \eqref{eq:1.17}, \eqref{eq:4.20}, and \eqref{eq:4.21} 
we obtain 
\begin{equation}
\label{eq:4.22}
\sup_{0<\tau<t}\|u_n(\tau)\|_{1,\rho}=\Psi_{1,\rho}[u_n](t)\le C\|\varphi\|_{1,\rho},
\quad
\|u_n(t)\|_{L^\infty(\Omega)}\le Ct^{-\frac{N}{2}}\|\varphi\|_{1,\rho},
\end{equation}
for all $0<t\le\mu\rho^2$ and $n=1,2,\dots$. 
Furthermore, applying Lemma~\ref{Lemma:4.3} with $v=u_m$ and $w=u_n$ 
and taking a sufficiently small $\mu $ if necessary, 
we see that 
$$
\sup_{0<\tau<\mu\rho^2}\|u_m-u_n\|_{1,\rho}\le 2M\|u_m(0)-u_n(0)\|_{1,\rho}. 
$$
Then, by the same argument as in the proof for the case $r>1$ 
we see that there exists a $L^1_{uloc}(\Omega)$-solution $u$ of \eqref{eq:1.1} 
in $\Omega\times[0,\mu\rho^2]$ satisfying \eqref{eq:1.9} and \eqref{eq:1.10}. 
Thus the proof of Theorem~\ref{Theorem:1.1} in the case $r=1$ is complete.
$\Box$\vspace{5pt}

\noindent
{\bf Proof of Theorem~\ref{Theorem:1.2} in the case $r=1$.} 
Let $v$ and $w$ be $L^1_{uloc}(\Omega)$-solutions of \eqref{eq:1.1} in $\Omega\times[0,T)$, where $0<T\le\infty$. 
Assume \eqref{eq:1.11}. 
Then, for any $0<T'<T$, 
we have 
$$
\|v(t)\|_{L^\infty(\Omega)}+\|w(t)\|_{L^\infty(\Omega)}\le Ct^{-\frac{1}{2(p-1)}},\qquad 0<t\le T'. 
$$
By Lemma~\ref{Lemma:4.3} 
we can find a positive constant $\mu\in(0,1)$ such that 
$$
\|(v(t)-w(t))_+\|_{1,\rho}\le 2M\|(v(0)-w(0))_+\|_{1,\rho}=0
$$
for all $0<t\le\min\{T',\mu\rho^2\}$. 
Repeating this argument, we see that 
$$
\|(v(t)-w(t))_+\|_{1,\rho}\le 0
$$
for all $0<t\le T'$. Since $T'$ is arbitrary, we deduce that $v(x,t)\le w(x,t)$ in $\Omega\times(0,T)$. 
Thus Theorem~\ref{Theorem:1.2} in the case $r=1$ follows. 
$\Box$\vspace{5pt}

\noindent
{\bf Proof of Corollary~\ref{Corollary:1.1}.} 
Let $p>1+1/N$ and $\varphi\in L^{N(p-1)}(\Omega)$. 
By \eqref{eq:2.2} we can find $\rho\in(0,\rho_*)$ such that 
$$
\|\varphi\|_{N(p-1),\rho}\le\gamma_1,
$$
where $\gamma_1$ is the constant given in Theorem~\ref{Theorem:1.1}. 
Then assertion~(i) follows from Theorem~\ref{Theorem:1.1}. 
Furthermore, if $\rho_*=\infty$ and $\varphi$ satisfies \eqref{eq:1.13}, 
then assertion~(i) of Theorem~\ref{Theorem:1.1} holds for any $\rho>0$. 
This implies assertion~(ii), and Corollary~\ref{Corollary:1.1} follows. 
$\Box$
\section{Applications}
In this section, as an application of Theorem~\ref{Theorem:1.1}, 
we give lower estimates of the blow-up time and  the blow-up rate for problem~\eqref{eq:1.1}.  
\subsection{Blow-up time}
%
Let $T(\lambda\psi)$ be the blow-up time of the solution of \eqref{eq:1.1} 
with the initial function $\varphi=\lambda\psi$. 
In this subsection we study the behavior of $T(\lambda\psi)$ as $\lambda\to\infty$ or $\lambda\to 0$. 
\begin{theorem}
\label{Theorem:5.1}
Let  $N\ge 1$ and $\Omega\subset{\bf R}^N$ be a uniformly regular domain of class $C^1$. 
Let $r$ satisfy 
$$
N(p-1)<r\le\infty\quad\mbox{if}\quad p\ge p_*
\qquad\mbox{and}\qquad
1\le r\le\infty\quad\mbox{if}\quad 1<p<p_*. 
$$
Then, for any $\psi\in L^r_{uloc,\rho}(\Omega)$ with $\rho>0$, 
there exists a positive constant $C$ such that 
$$
T(\lambda\psi)\ge
\left\{
\begin{array}{ll}
C(\lambda\|\psi\|_{r,\rho})^{-\frac{2r(p-1)}{r-N(p-1)}} & \mbox{if}\quad r<\infty,\vspace{3pt}\\
C(\lambda\|\psi\|_{L^\infty(\Omega)})^{-2(p-1)} & \mbox{if}\quad r=\infty,
\end{array}
\right.
$$
for all sufficiently large $\lambda$. 
\end{theorem}
{\bf Proof.}
Let $\gamma_1$ and $\mu$ be constants given in Theorem~\ref{Theorem:1.1}. 
If $r<\infty$, by Theorem~\ref{Theorem:1.1} we see that 
$$
T(\lambda\psi)\ge\mu\left(\frac{\gamma_1}{\lambda\|\psi\|_{r,\rho}}\right)^{2(\frac{1}{p-1}-\frac{N}{r})^{-1}}
\ge C(\lambda\|\psi\|_{r,\rho})^{-\frac{2r(p-1)}{r-N(p-1)}}
$$
for all sufficiently large $\lambda$. 
If $r=\infty$, then 
$$
\|\lambda\psi\|_{N(p-1),\rho}\le C\lambda\|\psi\|_{L^\infty(\Omega)}\rho^{\frac{1}{p-1}}. 
$$
It follows from Theorem~\ref{Theorem:1.1} that 
$$
T(\lambda\psi)\ge\mu\left(\frac{\gamma_1}{C\lambda\|\psi\|_{L^\infty(\Omega)}}\right)^{2(p-1)}
\ge C(\lambda\|\psi\|_{L^\infty(\Omega)})^{-2(p-1)}
$$
for all sufficiently large $\lambda$. 
Thus Theorem~\ref{Theorem:5.1} follows.
$\Box$
\begin{theorem}
\label{Theorem:5.2}
Let  $N\ge 1$ and $\Omega\subset{\bf R}^N$ be a uniformly regular domain of class $C^1$. 
Assume 
\begin{equation}
\label{eq:5.1}
\sup_{x\in\overline{\Omega}}|x|^\beta|\psi(x)|<\infty,
\end{equation}
where 
$0\le\beta<N$ if $1<p<p_*$ and $0\le\beta<1/(p-1)$ if $p\ge p_*$. 
Then there exists a positive constant $C_1$ such that 
\begin{equation}
\label{eq:5.2}
T(\lambda\psi)\ge C_1\lambda^{-\frac{2(p-1)}{1-\beta(p-1)}}
\end{equation}
for all sufficiently large $\lambda$. 
Furthermore, if $\Omega={\bf R}^N_+$ and 
\begin{equation}
\label{eq:5.3}
\inf_{x\in\Omega(0,\delta)}|x|^\beta\psi(x)>0
\end{equation}
for some $\delta>0$, then there exists a positive constant $C_2$ such that
\begin{equation}
\label{eq:5.4}
T(\lambda\psi)\le C_2\lambda^{-\frac{2(p-1)}{1-\beta(p-1)}}
\end{equation}
for all sufficiently large $\lambda$. 
\end{theorem}
{\bf Proof.} 
In the case $1<p<p_*$, 
let $r>1$, $r>N(p-1)$ and $\beta<N/r$. 
In the case $1<p<p_*$, let $r=1$. 
It follows from \eqref{eq:5.1} that 
$\rho^{\frac{1}{p-1}-\frac{N}{r}}\|\psi\|_{r,\rho}\le C\rho^{\frac{1}{p-1}-\beta}$
for all sufficiently small $\rho>0$. 
This together with Theorem~\ref{Theorem:1.1} implies \eqref{eq:5.2}. 

Assume \eqref{eq:5.3}. Let $v$ be a solution of
$$
\left\{
\begin{array}{ll}
\partial_t v=\Delta v & \quad\mbox{in}\quad{\bf R}^N_+\times(0,\infty),\vspace{3pt}\\
\nabla v\cdot\nu(x)=0 & \quad\mbox{in}\quad\partial{\bf R}^N_+\times(0,\infty),\vspace{3pt}\\
v(x,0)=A|x|^{-\beta}\chi_{B(0,\delta)} & \quad\mbox{in}\quad{\bf R}^N_+,
\end{array}
\right.
$$
where $A$ is a positive constant to be chosen as $\psi(x)\ge v(x,0)$ in ${\bf R}^N_+$. 
By \cite[Lemma~2.1.2]{DFL} we can find a constant $c_p$ depending only on $p$ such that 
\begin{equation}
\label{eq:5.5}
\lambda\|v(\cdot,0,t)\|_{L^\infty({\bf R}^{N-1})}
\le c_pt^{-\frac{1}{2(p-1)}},\qquad 0<t<T(\lambda v(0)). 
\end{equation}
On the other hand, since $T(\lambda\psi)\le T(\lambda v(0))$ and 
$$
\|v(\cdot,0,t)\|_{L^\infty({\bf R}^{N-1})}\ge Ct^{-\frac{\beta}{2}},\quad 0<t\le1,
$$
we have 
$$
\lambda T(\lambda\psi)^{\frac{1}{2(p-1)}-\frac{\beta}{2}}
\le \lambda T(\lambda v(0))^{\frac{1}{2(p-1)}-\frac{\beta}{2}}\le Cc_p,
$$
which implies \eqref{eq:5.4}. Thus Theorem~\ref{Theorem:5.2} follows.
$\Box$
\begin{remark}
For the case $\Omega=(0,\infty)$,  
Fern\'andez Bonder and Rossi {\rm\cite{FR}} proved 
$$
\lim_{\lambda\to\infty}\lambda^{2(p-1)}T(\lambda\psi)=T(\psi(0))
$$ 
provided that $\psi$ is bounded continuous and positive on $[0,\infty)$. 
\end{remark}
Motivated by \cite{LN},
we consider the case $\Omega={\bf R}^N_+$ and 
study the behavior of the blow-up time $T(\lambda\psi)$ as $\lambda\to 0$.
\begin{theorem}
\label{Theorem:5.3} 
Let $\Omega={\bf R}^N_+$ and assume 
\begin{equation}
\label{eq:5.6}
\sup_{x\in{\bf R}^N_+}\,(1+|x|)^\beta|\psi(x)|<\infty
\end{equation}
for some $\beta\ge 0$.  
Let $\lambda>0$ and consider problem~\eqref{eq:1.1} with $\varphi=\lambda\psi$. 
Then there exists a positive constant $C_1$ such that 
\begin{equation}
\label{eq:5.7}
T(\lambda\psi)\ge C_1f(\lambda)
\end{equation}
for all sufficiently small $\lambda>0$, where 
$$
f(\lambda):=\left\{
\begin{array}{lll}
\lambda^{-\frac{2(p-1)}{1-\beta(p-1)}} & \mbox{if}\quad p\ge p_*, & 0\le\beta<\frac{1}{p-1},\\
\lambda^{-\frac{2(p-1)}{1-\beta(p-1)}}  & \mbox{if}\quad 1<p<p_*, & 0\le\beta<N,\\
(\lambda|\log\lambda|)^{-\frac{2(p-1)}{1-N(p-1)}} & \mbox{if}\quad 1<p<p_*, & \beta=N,\\
\lambda^{-\frac{2(p-1)}{1-N(p-1)}} & \mbox{if}\quad1<p<p_*, & \beta>N.
\end{array}
\right. 
$$
Furthermore, if 
$$
\inf_{x\in{\bf R}^N_+}\,(1+|x|)^\beta\psi(x)>0,
$$
then there exists a positive constant $C_2$ such that 
\begin{equation}
\label{eq:5.8}
T(\lambda\psi)\le C_2f(\lambda)
\end{equation}
for all sufficiently small $\lambda>0$. 
\end{theorem}
{\bf Proof.}
Consider the case $p\ge p_*$. 
Let $0\le \beta<1/(p-1)$, $r>N(p-1)$ and $\beta<N/r$. 
By~\eqref{eq:5.6} we have 
$$
\rho^{\frac{1}{p-1}-\frac{N}{r}}\|\lambda\psi\|_{r,\rho}
\le C\lambda\rho^{-\beta+\frac{1}{p-1}}
$$
for all sufficiently large $\rho$. 
Similarly, in the case $p<p_*$, 
it follows from \eqref{eq:5.6} that 
$$
\rho^{\frac{1}{p-1}-N}\|\lambda\psi\|_{1,\rho}
\le 
\left\{
\begin{array}{ll}
 C\lambda\rho^{-\beta+\frac{1}{p-1}}
 & \mbox{if}\quad 0\le\beta<N,\vspace{3pt}\\
 C\lambda\rho^{\frac{1}{p-1}-N}\log\rho & \mbox{if}\quad \beta=N,\vspace{3pt}\\
 C\lambda\rho^{\frac{1}{p-1}-N}
 & \mbox{if}\quad\beta>N,
\end{array}
\right.
$$
for all sufficiently large $\rho$. 
Therefore, by Theorem~\ref{Theorem:1.1} 
we obtain in the case $p\ge p_*$ 
$$
T(\lambda\psi)\ge C\lambda^{-\frac{2}{-\beta+\frac{1}{p-1}}}
=C\lambda^{-\frac{2(p-1)}{1-\beta(p-1)}}
$$
and in the case $1<p<p_*$ 
$$
T(\lambda\psi)
\ge\left\{
\begin{array}{ll}
C\lambda^{-\frac{2(p-1)}{1-\beta(p-1)}} & \mbox{if}\quad 0\le\beta<N,\vspace{3pt}\\
C(\lambda|\log\lambda|)^{-\frac{2(p-1)}{1-N(p-1)}} & \mbox{if}\quad\beta=N,\vspace{3pt}\\
C\lambda^{-\frac{2(p-1)}{1-N(p-1)}} & \mbox{if}\quad\beta>N,
\end{array}
\right.
$$
for all sufficiently small $\lambda>0$. 
These imply \eqref{eq:5.7}. 

Let $v$ be a solution of 
$$
\left\{
\begin{array}{ll}
\partial_t v=\Delta v & \quad\mbox{in}\quad{\bf R}^N_+\times(0,\infty),\vspace{3pt}\\
\nabla v\cdot\nu(x)=0 & \quad\mbox{in}\quad\partial{\bf R}^N_+\times(0,\infty),\vspace{3pt}\\
v(x,0)=A(1+|x|)^{-\beta} & \quad\mbox{in}\quad{\bf R}^N_+,
\end{array}
\right.
$$
where $A$ is a positive constant to be chosen as $\psi(x)\ge v(x,0)$ in ${\bf R}^N_+$. 
Since $T(\lambda\psi)\le T(\lambda v(0))$ and 
$$
\|v(\cdot,0,t)\|_{L^\infty({\bf R}^{N-1})}
\ge\left\{
\begin{array}{ll}
Ct^{-\frac{\beta}{2}} & \mbox{if}\quad 0\le\beta<N,\vspace{3pt}\\
Ct^{-\frac{N}{2}}\log t & \mbox{if}\quad\beta=N,\vspace{3pt}\\
Ct^{-\frac{N}{2}} & \mbox{if}\quad\beta>N,
\end{array}
\right.
$$
for all sufficiently large $t$, 
by a similar argument as in the proof of \eqref{eq:5.4} we obtain \eqref{eq:5.8}. 
Thus Theorem~\ref{Theorem:5.3} follows.
$\Box$
\subsection{Blow-up rate}
Let $u$ be a solution of \eqref{eq:1.1} in $\Omega\times[0,T)$, where $0<T<\infty$, 
such that $u$ blows up at $t=T$. 
In this subsection, as a corollary of Theorem~\ref{Theorem:1.1}, 
we state a result on lower estimates of the blow-up rate of the solution $u$. 
Blow-up rate of positive solutions for problem~\eqref{eq:1.1} was first obtained by Fila and Quittner \cite{FQ}, 
where it was shown that 
\begin{equation}
\label{eq:5.9}
\limsup_{t\to T}\,(T-t)^{\frac{1}{2(p-1)}}\|u(t)\|_{L^\infty(\Omega)}<\infty
\end{equation}
holds in the case where $\Omega$ is a ball, the initial function $\varphi$ is radially symmetric and 
satisfies some monotonicity assumptions. 
Subsequently, it was proved that \eqref{eq:5.9} holds for positive solutions in the following cases: 
\begin{itemize}
  \item $\Omega$ is a bounded smooth domain, $(N-2)p<N$ and $\partial_t u\ge 0$ in $\Omega\times(0,T)$ 
  (see \cite{GH}, \cite{H1} and \cite{HM});
  \item $\Omega$ is a bounded smooth domain and $p\le 1+1/N$ (see \cite{H3});
  \item $\Omega={\bf R}^N_+$ and $(N-2)p<N$ (see \cite{CF02}). 
\end{itemize}
See \cite{QS2} for sign changing solutions. 
On the other hand, for positive solutions, 
it was shown in \cite{HM} that 
\begin{equation}
\label{eq:5.10}
\liminf_{t\to T}\,(T-t)^{\frac{1}{2(p-1)}}\|u(t)\|_{L^\infty(\Omega)}>0 
\end{equation}
holds if $\Omega$ is a bounded smooth domain
(see also \cite{GH} and \cite{H1}). 
\vspace{5pt}

We state a result on lower estimates of the blow-up rate of the solutions. 
Theorem~\ref{Theorem:5.4} is a generalization of \eqref{eq:5.10} 
and it holds without the boundedness of the domain $\Omega$ and the positivity of the solutions. 
\begin{theorem}
\label{Theorem:5.4} 
Let $N\ge 1$ and $\Omega\subset{\bf R}^N$ be a uniformly regular domain of class $C^1$. 
Let $u$ be a solution of \eqref{eq:1.1} blowing up at $t=T<\infty$. 
Then 
\begin{equation}
\label{eq:5.11}
\liminf_{t\to T}\,(T-t)^{\frac{1}{2(p-1)}-\frac{N}{2r}}\|u(t)\|_{L^r(\Omega)}>0,
\end{equation}
where 
\begin{equation}
\label{eq:5.12}
\left\{
\begin{array}{ll}
N(p-1)\le r\le\infty & \mbox{if}\quad p>1+1/N,\vspace{3pt}\\
1<r\le\infty & \mbox{if}\quad p=1+1/N,\vspace{3pt}\\
1\le r\le\infty & \mbox{if}\quad 1<p<1+1/N. 
\end{array}
\right.
\end{equation}
\end{theorem}
{\bf Proof.}
Let $1\le r<\infty$ satisfy \eqref{eq:5.12}. 
By Theorem~\ref{Theorem:1.1} we can find positive constants $\gamma_1$ and  $\mu$ such that, if 
$$
||u(T-t)||_{r,\rho}\le\gamma_1\rho^{\frac{N}{r}-\frac{1}{p-1}}
$$
for some $\rho\in(0,\rho_*/2)$, 
then the solution $u$ exists in $\Omega\times(0,T-t+\mu\rho^2]$. 
Since the solution $u$ blows up at $t=T$, 
we can find a constant $\delta>0$ such that 
\begin{equation}
\label{eq:5.13}
||u(T-t)||_{r,\rho(t)}>\gamma_1\rho(t)^{\frac{N}{r}-\frac{1}{p-1}}\quad\mbox{for}\quad t\in(T-\delta,T),
\end{equation}
where
$$
\rho(t):=\left(\frac{T-t}{\mu}\right)^{\frac{1}{2}}. 
$$
This implies \eqref{eq:5.11} in the case $r<\infty$. 
Furthermore, by \eqref{eq:5.13}, 
for any $t\in(T-\delta, T)$, 
there exist $x(t)\in\overline{\Omega}$ and $y(t)\in\Omega(x(t),\rho(t))$ such that 
$$
C\rho(t)^Nu(y(t),t)^r\ge \int_{\Omega(x(t),\rho(t))}u(y,t)^r\,dy\ge\frac{\gamma_1}{2}\rho(t)^{N-\frac{r}{p-1}}. 
$$
This yields \eqref{eq:5.11} in the case $r=\infty$, and 
Theorem~\ref{Theorem:5.4} follows. 
$\Box$
\vspace{5pt}

\noindent
{\bf Acknowledgements.}
The first author was supported in part 
by the Grant-in-Aid for for Scientific Research (B)(No.~23340035),
from Japan Society for the Promotion of Science.
The second author was supported in part by Research Fellow 
of Japan Society for the Promotion of Science.
\bibliographystyle{amsplain}

\end{document}